\documentclass[preprint]{elsarticle}
\usepackage{lineno,hyperref}
\modulolinenumbers[5]

\journal{Elsevier}

\usepackage{amsmath}
\usepackage{amssymb}
\usepackage{amsthm}
\usepackage{amscd}
\usepackage{amsfonts}
\usepackage{graphicx}%

\usepackage{bm}

\usepackage{float}

\usepackage{color}
\usepackage{setspace}
\usepackage{slashbox}

\theoremstyle{plain} 
\numberwithin{equation}{section}
\newtheorem{theorem}{Theorem}[section]

\theoremstyle{definition}

\newtheorem{remark}[theorem]{Remark}

\def\x{\mathbf{x}}
\def\y{\mathbf{y}}
\def\D{\mathcal{D}}
\def\div{\mathrm{div\,}}
\def\q{\mathbf{q}}
\def\p{\mathbf{p}}

\def\w{\mathbf{w}}
\def\E{\mathbb{E}}

\def\r{\mathbf{p}}

\begin{document}
\begin{frontmatter}

\title{Energy Conserving Galerkin Approximation of Two Dimensional Wave Equations with Random Coefficients}

\author[myaddress]{Ching-Shan Chou \fnref{fn1}}

\author[myaddress]{Yukun Li}

\author[myaddress]{Dongbin Xiu}

\address[myaddress]{Department of Mathematics, The Ohio State University }

\fntext[fn1]{CSC is supported by NSF grant DMS-1253481}

\begin{abstract}
Wave propagation problems for heterogeneous media are known to have many applications in physics and engineering. 
Recently, there has been an increasing interest in stochastic effects due to the uncertainty,
which may arise from impurities of the media. This work considers a two-dimensional wave equation 
with random coefficients which may be discontinuous in space. 
Generalized polynomial chaos method is used in conjunction with stochastic Galerkin approximation, and
local discontinuous Galerkin method is used for spatial discretization.
Our method is shown to be energy preserving in semi-discrete form as well as in
fully discrete form, when leap-frog time discretization is used. 
Its convergence rate is proved to be optimal and the error grows linearly in time. 
The theoretical properties of the proposed scheme are validated by numerical tests. 
 \end{abstract}

\end{frontmatter}


\noindent{\bf{AMS}}
65N12, 
65N15, 
65N30. 

\linenumbers

\section{Introduction}\label{sec-1}

Consider the following second order deterministic wave equations
\begin{align*}
\frac{\partial^2 u(t,\x)}{\partial t^2} &=  \mathrm{div}(a^2(\x)\nabla u(t,\x)),\qquad \x\in\D,\ t\in\mathcal{T},\\
u(\x,0)&=u_0(\x),\qquad u_t(\x,0)=v_0(\x),
\end{align*}
subject to homogeneous Dirichlet or periodic boundary conditions. Here $\D$ denotes a two-dimensional physical domain, $\mathcal{T}$ denotes a time range, and $a(\x)$ denotes the speed of wave propagation. 
An important property of the wave equation is its conservation of energy. Therefore, recently there is an increasing interest in energy conserving numerical methods for wave equations, and it has been shown that these methods preserve the shape and phase of smooth shaped waves. 

Here we focus on discontinuous Galerkin (DG) method for discretization in physical space. Historically, there are basically two approaches to design energy conserving DG methods. One approach is to use staggered meshes. Chung and Engquist have used this approach and proposed an optimal and energy conserving DG scheme for the first-order wave equation \cite{chung2006optimal,chung2009optimal}. 
The other approach is to use the central numerical flux in DG method \cite{fezoui2005convergence}. However, the convergence for this scheme is suboptimal theoretically, and numerically shown to be optimal/suboptimal for even/odd degree polynomial basis \cite{fezoui2005convergence}.
As an alternative, Xing and Chou developed a local discontinuous Galerkin (LDG) (\cite{chou2014optimal, xing2013energy}) that produces both energy conservation and optimal convergence rate.

In practical applications, the wave propagation speed $a$ is unlikely to be deterministic, because the media in which the wave propagates often 
have random impurities. This leads us to consider $a$ as a function of both space and random variables, and its associated solution $u$, a function of space, time and random variables. To characterize the stochastic function $u$, a popular and robust approach is Monte-Carlo method. As a brute-force sample-based method, a large number of samples are usually needed to achieve satisfactory accuracy, and therefore it is known to be  computationally expensive. One efficient alternative is polynomial chaos (PC) approximation, originally developed by Ghanem and Spanos using Wiener-Hermite expansion and finite element discretization for a range of problems \cite{ghanem2003stochastic}. It was later extended by Xiu and Karniadakis \cite{xiu2002wiener} to generalized polynomial chaos (gPC) expansion, in which general orthogonal polynomials were considered. Based on gPC expansion and stochastic Galerkin projection, the original random PDE can be transformed into a system of deterministic equations which can be solved by existing numerical methods \cite{babuska2004galerkin,ghanem2003stochastic,frauenfelder2005finite,xiu2002wiener}. 
Among the existing work, the stochastic Galerkin methods for the first-order random hyperbolic problems were considered in \cite{gottlieb2008galerkin,jin2016well, wu2017stochastic}. On a different front, stochastic collocation methods have also been considered for scalar hyperbolic equations (\cite{tangzhou2010}) and second-order wave equation with a discontinuous random speed (\cite{motamed2014}). Stochastic Galerkin and stochastic collocation are the two main approaches for
problems with random inputs. They have different properties and both are useful for different problems. Their comparison is beyond the scope of this paper. Here we focus
on the properties of stochastic Galerkin method for wave equations particularly in conjunction with LDG method for energy conservation.

In this paper, we apply the gPC Galerkin framework, along with LDG,  to the second-order wave equation directly, without transforming it into a first order hyperbolic system.  Our method is thus a Galerkin approximation in both physical space and random space.  More importantly, we demonstrate that the resulting numerical scheme is energy conserving.
Consequently, it induces much less errors for long time integration. 
%
We first examine the stability of the stochastic wave equation, with respect to the random wave speed $a$ by characterizing its solution dependence on the random coefficient. This is similar to the previous work for the elliptic problem \cite{lord2014introduction}. 
Upon presenting the detail of the numerical scheme,
we then prove that the numerical scheme is energy conserving  in both semi-discrete and fully discrete forms. 
Finally, we show that by
taking a suitable projection for the initial conditions, our numerical scheme achieves optimal convergence rate.

The paper is organized as follows. In Section 2, the stability of the problem with respect to the random coefficient $a$ is proved.
In Section 3, we present our numerical method of gPC expansion and LDG framework. The energy conserving properties are proved for
both semi-discrete and fully-discrete (leap-frog) schemes.
In Section 4, error estimates are presented for the semi-discrete numerical method.
In Section 5, we present numerical tests with random $a$, continuous or discontinuous in space, to demonstrate the 
energy conserving properties and error estimates proved in previous sections.
Concluding remarks are given in Section 6.

\section{Dependence of Solution on Random Wave Speed}\label{sec-2}

In this paper, consider the following two-dimensional wave equation with random coefficient 
\begin{equation}\label{eq1}
\frac{\partial^2 u(t,\x,\y)}{\partial t^2} =  \mathrm{div}(a^2(\x,\y)\nabla u(t,\x,\y)), 
\end{equation}
where $\x$ denotes the spatial variables in the two-dimensional domain $\mathcal{D}$ and $\y=(y_1,y_2,\ldots,y_N)\in \mathbb{R}^N,\ N\ge1$, is a random vector with independent and identically distributed components. 
Equation (\ref{eq1}) is subject to initial condition
\begin{equation}\label{eq2}
u(0,\x,\y) = u_0(\x,\y),\qquad u_t(0,\x,\y) = v_0(\x,\y),
\end{equation}
and the homogeneous Dirichlet boundary conditions
\begin{equation}\label{eq3}
u(t,\x,\y) = 0\qquad \x\in\partial\mathcal{D}.
\end{equation}
The coefficient $a^2(\x,\y)$ is assumed to be positive for all $\x$ and $\y$. 
Because $a(\x,\y)$ is associated with the media in which the wave propagates, 
Eq.~\eqref{eq1} models wave propagation in heterogeneous media subject to random variations.
For the convenience of applying the LDG framework later, we first rewrite \eqref{eq1} into the equivalent system
\begin{alignat}{2}\label{eq_u}
\frac{\partial^2 u(t,\x,\y)}{\partial t^2} &=  \mathrm{div}(a(\x,\y)\q(t,\x,\y)),\qquad &&\\
\q(t,\x,\y)&=a(\x,\y)\nabla u(t,\x,\y)\qquad &\q\in \mathbb{R}^{2\times 1}.&\label{eq_q}
\end{alignat}

In this section, we would like to establish the stability of Eqs.~\eqref{eq_u} and \eqref{eq_q} with respect to the wave speed
coefficient $a(\x,\y)$; in other words, we will show that if a small perturbation is made on $a$, in either $\x$ or $\y$, the solution will be close
to that without perturbation.
The stability of the problem is relevant because in real applications, the function $a(\x,\y)$ may be approximated and not exact. Hence it is necessary 
to show that as long as the approximation on $a$ is sufficiently accurate, the resulting solution will be sufficiently close to the exact solution.



First, we take the time derivative of \eqref{eq_q},
\begin{align}
\q_t(t,\x,\y)&=  a(\x,\y)\nabla u_{t}(t,\x,\y).\label{eq_qt}
\end{align}
After taking the expectation with respect to $\y$ on both sides of the weak form of \eqref{eq_u} and \eqref{eq_qt}, we obtain the following: $u\in L^{2}(\mathbb{R}^N;{H}^2(\mathcal{T};{H}^{-1}( \mathcal{D})))\cap L^{2}(\mathbb{R}^N;{L}^2(\mathcal{T};{H}^1_0(\mathcal{D})))$ and $\q\in (\mathbf{L}^{2}(\mathbb{R}^N;\mathbf{L}^2(\mathcal{T};\mathbf{H}^1( \mathcal{D}))))^2$ satisfy
\begin{align}\label{eq40}
\mathbb{E} \left[(u_{tt},p)_\mathcal{D}\right] +\mathbb{E}\left[(a\q,\nabla p)_\mathcal{D}\right] &= 0\qquad \forall {p}\in L^{2}(\mathbb{R}^N;{H}^1_0( \mathcal{D})),\\
\mathbb{E}\left[(\q_t,\w)_\mathcal{D}\right] +  \mathbb{E}\left[(au_t,\div\w )_\mathcal{D}\right]+  \mathbb{E}\left[(\nabla au_t,{\w})_\mathcal{D}\right]&= 0\qquad \forall {\w}\in (\mathbf{L}^{2}(\mathbb{R}^N;\mathbf{H}^1( \mathcal{D})))^2,\label{eq41}
\end{align}
where ${H}^1_0( \mathcal{D})$ denotes the set of functions in ${H}^1( \mathcal{D})$ with vanishing boundary values. Here we use $(\cdot\ ,\ \cdot)_\mathcal{D}$ to denote the integral of the product (inner product) over $\mathcal{D}$ if the arguments are scalar (vector) functions.

Suppose $\widetilde{a}(\x,\y)$ is a perturbed function of $a(\x,\y)$, and its corresponding solutions are $\widetilde{u}(t,\x,\y)$ and $\widetilde{\q}(t,\x,\y)$. Then $\widetilde{u}$ and $\widetilde{\q}$ satisfy
\begin{align}\label{eq43}
\mathbb{E} \left[(\widetilde{u}_{tt},p)_\mathcal{D}\right] +\mathbb{E}\left[(\widetilde{a}\widetilde{\q},\nabla p)_\mathcal{D}\right] &= 0\qquad \forall {p}\in L^{2}(\mathbb{R}^N;{H}_0^1( \mathcal{D})),\\
\mathbb{E}(\widetilde{\q}_t,\w) +  \mathbb{E}(\widetilde{a}\widetilde{u}_t,\div\w ) +  \mathbb{E}(\nabla \widetilde{a}(\widetilde{u})_t,{\w})&= 0\qquad \forall {\w}\in (\mathbf{L}^{2}(\mathbb{R}^N;\mathbf{H}^1( \mathcal{D})))^2.\label{eq44}
\end{align}
%
%
We assume that both $a$ and $\widetilde{a}$ are bounded from above and from
below away from 0, that is,
$a^2(\x,\y)$ and $\widetilde{a}^2(\x,\y)$ belong to $L^{\infty}(\mathbb{R}^N,W^{1,\infty}(\mathcal{D}))$ and 
\begin{align}
0<a_{\min}\le \|a^2(\x,\y)\|_{L^{\infty}(\mathbb{R}^N;W^{1,\infty}(\mathcal{D}))} \le a_{\max}<+\infty\qquad \text{a.e. in}\ \mathcal{D}\times\mathbb{R}^N,\notag\\
0<\widetilde{a}_{\min}\le \|\widetilde{a}^2(\x,\y)\|_{L^{\infty}(\mathbb{R}^N;W^{1,\infty}(\mathcal{D}))} \le \widetilde{a}_{\max}<+\infty\qquad \text{a.e. in}\ \mathcal{D}\times\mathbb{R}^N.\notag
\end{align}
Based on the above assumptions, and assuming that $a(\x,\y)$ and $\widetilde{a}(\x,\y)$ have the same sign, we can easily show that given an arbitrary $\epsilon>0$, if
\begin{align}\label{20161206_add1}
\|a^2(\x,\y)-\widetilde{a}^2(\x,\y)\|_{L^{\infty}(\mathbb{R}^N;W^{1,\infty}(\mathcal{D}))}\le\epsilon,
\end{align}
then 
\begin{align}\notag
\|a(\x,\y)-\widetilde{a}(\x,\y)\|_{L^{\infty}(\mathbb{R}^N;W^{1,\infty}(\mathcal{D}))}\le C_1\epsilon,
\end{align}
where $C_1 = C(\sqrt{\widetilde{a}_{\max}}+\sqrt{a_{\max}})\slash(\sqrt{\widetilde{a}_{\min}}+\sqrt{a_{\min}})^2$.

We define the difference between the solutions of the perturbed and the original systems to be $\delta_u=u-\widetilde{u}$ and $\delta_\q=\q-\widetilde{\q}$.
In the following theorem, we prove the bound of the averaged $L^2$ norm of the difference between the solutions in terms of the perturbation in the coefficient $a(\x,\y)$.

\begin{theorem}\label{thm20170318_1}
Let $u(t,\x,\y)$ and $\widetilde{u}(t,\x,\y)$ be solutions of \eqref{eq40}--\eqref{eq41} and \eqref{eq43}--\eqref{eq44}, respectively.
If the initial conditions satisfy 
\begin{align*}
\bigl(\mathbb{E}[\|u_t(0,\x,\y)- \widetilde{u}_t(0,\x,\y)\|_{L^2(\mathcal{D})}^2] + \mathbb{E}[\|\q(0,\x,\y)-\widetilde{\q}(0,\x,\y)\|_{L^2(\mathcal{D})}^2]\bigr)^{\frac12} \le C\epsilon,
\end{align*}
then we have 
\begin{align}\notag
\bigl(\mathbb{E}[\|(\delta_u)_{t}\|_{L^2(\mathcal{D})}^2]\bigr)^\frac12 + \bigl(\mathbb{E}[\|\delta_\q\|_{L^2(\mathcal{D})}^2]\bigr)^{\frac12} \le C(t+1)\epsilon.
\end{align}
\end{theorem}

\begin{proof}
Subtracting \eqref{eq43}--\eqref{eq44} from \eqref{eq40}--\eqref{eq41} respectively, we have
\begin{align}\label{eq45}
\mathbb{E}\left[ ((\delta_u)_{tt},p)_{\mathcal{D}}\right]+\mathbb{E}\left[(a\q-\widetilde{a}\widetilde{\q},\nabla p)_{\mathcal{D}}\right] &= 0\qquad \forall {p}\in L^{2}(\mathbb{R}^N;{H}_0^1( \mathcal{D})),&\\
\mathbb{E}\left[((\delta_\q)_{t},\w)_\mathcal{D}\right] +  \mathbb{E}\left[(au_t-\widetilde{a}\widetilde{u}_t,\div\w )_\mathcal{D}\right]
+  \mathbb{E}\left[(\nabla au_t-\nabla \widetilde{a}\widetilde{u},{\w})_\mathcal{D}\right]&= 0\qquad \forall {\w}\in (\mathbf{L}^{2}(\mathbb{R}^N;\mathbf{H}^1( \mathcal{D})))^2.\label{eq46}
\end{align}
Choosing $p=(\delta_u)_t$ in \eqref{eq45}, $\w=\delta_\q$ in \eqref{eq46} and applying integration by parts to the second term of \eqref{eq45} yields
\begin{align}
\mathbb{E}\left[ ((\delta_u)_{tt},(\delta_u)_{t})_\mathcal{D}\right]+ \mathbb{E}\left[((\delta_\q)_t,\delta_\q)_\mathcal{D}\right]  +  \mathbb{E}\left[ (au_t-\widetilde{a}\widetilde{u}_t,\div \delta_\q )_\mathcal{D}\right]
+\mathbb{E}\left[ ((\nabla au_t-(\nabla \widetilde{a}\widetilde{u}_t,\delta_\q)_\mathcal{D} \right]&\notag\\
-\mathbb{E}\left[ (\nabla a\cdot\q-\nabla \widetilde{a}\cdot\widetilde{\q},(\delta_u)_t)_\mathcal{D}\right]-\mathbb{E}\left[(a\ \div\q-\widetilde{a}\ \div\widetilde{\q},(\delta_u)_t)_\mathcal{D} \right]&=0.\label{eq47}
\end{align}
Consider the fourth and the fifth terms on the left-hand side of \eqref{eq47}, we have
\begin{align}
&-\mathbb{E}\left[(\nabla au_t-\nabla \widetilde{a}\widetilde{u}_t,\delta_\q)_\mathcal{D}\right]
+\mathbb{E}\left[(\nabla a\cdot\q-\nabla \widetilde{a}\cdot\widetilde{\q},(\delta_u)_t)_\mathcal{D}\right]\notag\\
=&-\bigl(\mathbb{E}[(\nabla a\cdot \q,\widetilde{u}_t)_\mathcal{D}]-\mathbb{E}[(\nabla \widetilde{a}\cdot \q,\widetilde{u}_t)_\mathcal{D}\bigr)
-\bigl(\mathbb{E}[(\nabla \widetilde{a}\cdot \widetilde{\q},u_t)_\mathcal{D}]-\mathbb{E}[(\nabla a\cdot \widetilde{\q},u_t)_\mathcal{D}\bigr)\notag\\
=&\mathbb{E}[(\nabla(\widetilde{a}-a)\cdot \delta_\q,u_t)_\mathcal{D}]-\mathbb{E}[(\nabla(\widetilde{a}-a)\cdot \q,(\delta_u)_t)_\mathcal{D}]\notag\\
\le&\|\nabla(\widetilde{a}-a)\|_{L^{\infty}(\mathcal{D}\times\mathbb{R}^N)}\mathbb{E}[\|u_t\|_{L^{2}(\mathcal{D})}\|\delta_\q\|_{L^2(\mathcal{D})}]
+\|\nabla(\widetilde{a}-a)\|_{L^{\infty}(\mathcal{D}\times\mathbb{R}^N)}\mathbb{E}[\|\q\|_{L^{2}(\mathcal{D})}\|(\delta_u)_t\|_{L^2(\mathcal{D})}]\notag\\
\le&C\epsilon(\mathbb{E}[\|\delta_\q\|_{L^2(\mathcal{D})}^2]+\mathbb{E}[\|(\delta_u)_t\|_{L^2(\mathcal{D})}^2])^{1\slash2},\label{eq48}
\end{align}
where $C=C_1(\mathbb{E}[\|u_t\|^2_{L^{2}(\mathcal{D})}]+\mathbb{E}[\|\q\|^2_{L^{2}(\mathcal{D})})^{1\slash2}$.\\
Consider the third and the sixth term on the left-hand side of \eqref{eq47}, we have
\begin{align}
&-\mathbb{E}[(au_t-\widetilde{a}(\widetilde{u})_t,\div \delta_\q)_\mathcal{D}]
+\mathbb{E}[(a\div\q-\widetilde{a}\div\widetilde{\q},(\delta_u)_t)_\mathcal{D}]\notag\\
=&\mathbb{E}[(au_t(t,\x,\y),\div \widetilde{\q})_\mathcal{D}]+\mathbb{E}[(\widetilde{a}(\widetilde{u})_t,\div \q)_\mathcal{D}]
-\mathbb{E}[(\widetilde{a}u_t,\div \widetilde{\q})_\mathcal{D}]-\mathbb{E}[(a(\widetilde{u})_t,\div \q)_\mathcal{D}]\notag\\
=&-\mathbb{E}[(a-\widetilde{a},\div\q\,(\widetilde{u})_t-\div\widetilde{\q}\,u_t)_\mathcal{D}]\notag\\
=&-\mathbb{E}[(a-\widetilde{a},\div\q\,(\widetilde{u})_t-\div\q\,u_t)_\mathcal{D}]
-\mathbb{E}[(a-\widetilde{a},\div\q\,u_t-\div\widetilde{\q}\,u_t)_\mathcal{D}]\notag\\
\le&\|\widetilde{a}-a\|_{L^{\infty}(\mathcal{D}\times\mathbb{R}^N)}\mathbb{E}\bigl[\|\div\q\|_{L^{2}(\mathcal{D})}\|(\delta_u)_t\|_{L^2(\mathcal{D})}\bigr]\notag\\
&\quad+\|\widetilde{a}-a\|_{L^{\infty}(\mathcal{D}\times\mathbb{R}^N)}\mathbb{E}\bigl[\|\nabla u_t(t,\x,\y)\|_{L^{2}(\mathcal{D})}\|\delta_\q\|_{L^{2}(\mathcal{D})}\bigr]\notag\\
&\quad+\|\nabla(\widetilde{a}-a)\|_{L^{\infty}(\mathcal{D}\times\mathbb{R}^N)}\mathbb{E}\bigl[\|u_t(t,\x,\y)\|_{L^{2}(\mathcal{D})}\|\delta_\q\|_{L^{2}(\mathcal{D})}\bigr]\notag\\
\le&C\epsilon(\mathbb{E}[\|\delta_\q\|_{L^2(\mathcal{D})}^2]+\mathbb{E}[\|(\delta_u)_t\|_{L^2(\mathcal{D})}^2])^{1\slash2},\label{eq49}
\end{align}
where $C=C_1(\mathbb{E}[\|\nabla u_t\|^2_{L^{2}(\mathcal{D})}]+\mathbb{E}[\|u_t\|^2_{L^{2}(\mathcal{D})}]+\mathbb{E}[\|\div\q\|^2_{L^{2}(\mathcal{D})}])^{1\slash2}$.\\

By \eqref{eq47}-\eqref{eq49}, we have
\begin{align*}
\frac12\frac{d}{dt}\bigl(\mathbb{E}[\|(\delta_u)_{t}\|_{L^2(\mathcal{D})}^2] + \mathbb{E}[\|\delta_\q\|_{L^2(\mathcal{D})}^2]\bigr) \le C\bigl(\mathbb{E}[\|\delta_\q\|_{L^2(\mathcal{D})}^2]+\mathbb{E}[\|(\delta_u)_t\|_{L^2(\mathcal{D})}^2]\bigr)^{1\slash2}.
\end{align*}
and therefore
\begin{align*}
\frac{d}{dt}\bigl(\mathbb{E}[\|(\delta_u)_{t}\|_{L^2(\mathcal{D})}^2] + \mathbb{E}[\|\delta_\q\|_{L^2(\mathcal{D})}^2]\bigr)^{\frac12} \le C\epsilon.
\end{align*}
Because
\begin{align*}
\bigl(\mathbb{E}[\|(\delta_u)_t(0,\x,\y)\|_{L^2(\mathcal{D})}^2] + \mathbb{E}[\|\delta_\q(0,\x,\y)\|_{L^2(\mathcal{D})}^2]\bigr)^{\frac12} \le C\epsilon,
\end{align*}
we obtain
\begin{align}\notag
\bigl(\mathbb{E}[\|(\delta_u)_{t}\|_{L^2(\mathcal{D})}^2]\bigr)^\frac12 + \bigl(\mathbb{E}[\|\delta_\q\|_{L^2(\mathcal{D})}^2]\bigr)^{\frac12} \le C(t+1)\epsilon.
\end{align}
\end{proof}

\section{An Energy Conserving Numerical Method}\label{sec-2}
Assume that the solution of \eqref{eq_u}-\eqref{eq_q} can be expanded using polynomial chaos expansion
\begin{align}\label{eq54}
u(t,\x,\y)=\sum_{m=1}^{\infty}v_m(t,\x)\Phi_m(\y),\\
\q(t,\x,\y)=\sum_{m=1}^{\infty}\r_m(t,\x)\Phi_m(\y),\label{eq20170319_1}
\end{align}
where $\{\Phi_m(\y)\}_{m=1}^\infty$ are $N$-variate orthonormal polynomials, 
and the choice of the polynomials is based on the underlying probability density function $\rho(\y)$ for the random variable $\y$ \cite{xiu2002wiener}.
Specifically, 
\begin{equation}
\int \rho(\y)\Phi_{m}(\y)\Phi_{m'}(\y)d\y =\delta_{mm'},\label{eq8}
\end{equation}
where $\delta_{mm'}$ are the Kronecker delta functions. 
These orthonormal polynomials can be written as the products of univariate polynomials, 
\begin{equation}\label{eq6}
\Phi_{m}(\y) = \phi_{m_1}(y_1)\phi_{m_2}(y_2)\cdot\ldots\cdot\phi_{m_N}(y_N),
\end{equation}
with $m_i$ being the degree of $\phi(y_i)$ in the $y_i$-direction and $m$ the corresponding index integer for the vector index $(m_1,m_2,\cdots,m_N)$. 
$\rho(\y)$, the joint probability distribution function for $\y$, can be written as a product of univariate probability density function 
$\prod_{i=1}^N\rho_i(y_i)$, with $\rho_i(y_i)$ being the probability density function for $y_i$.

Substituting \eqref{eq54} and \eqref{eq20170319_1} into Eqs.~\eqref{eq_u}-\eqref{eq_q}, we have for all $k$
\begin{align}\label{eq20170319_2}
\frac{\partial^2 v_k}{\partial t^2}(t,\x)&=\sum_{j=1}^{\infty}\mathrm{div}(a_{kj}(\x) \r_j),\\
\r_k(t,\x)&=\sum_{j=1}^{\infty}a_{kj}(\x)\nabla v_j,\label{eq20170319_3}
\end{align}
where
\begin{align}\label{eq20170319_4}
a_{kj}(\x) = \int a(\x,\y)\Phi_k(\y)\Phi_j(\y)\rho(\y)d\y, \qquad j,k\ge 1.
\end{align}


If we look for the $P$-th order gPC approximation of $u$ and $\q$, i.e., 
\begin{align}\label{eq5}
u(t,\x,\y) \approx {u}_M(t,\x,\y):=\sum_{m=1}^M\widehat{v}_m(t,\x)\Phi_m(\y),\\
\q(t,\x,\y) \approx {\q}_M(t,\x,\y):=\sum_{m=1}^M\widehat{\r}_m(t,\x)\Phi_m(\y),\label{eq20160316_1}
\end{align}
where $M=\left(\begin{smallmatrix}N+P\\ N\end{smallmatrix}\right)$,
then by Galerkin projection, the coefficients in \eqref{eq5}-\eqref{eq20160316_1} satisfy
\begin{align}\label{eq9}
\frac{\partial^2 \widehat{v}_k}{\partial t^2}(t,\x)&=\sum_{j=1}^M\mathrm{div}(a_{kj}(\x) \widehat{\r}_j),\\
\widehat{\r}_k(t,\x)&=\sum_{j=1}^Ma_{kj}(\x)\nabla \widehat{v}_j,\label{eq20160316_4}
\end{align}
where $a_{kj}(\x)$ is defined in \eqref{eq20170319_4}.

We denote $\widehat{\mathbf{v}} = (\widehat{v}_1,\widehat{v}_2,\ldots,\widehat{v}_M)^T \in \mathbf{R}^{M\times 1}$ and $\widehat{\mathbf{S}} = (\widehat{\r}_1^T,\widehat{\r}_2^T,\ldots,\widehat{\r}_M^T)^T\in \mathbf{R}^{M\times 2}$. 
By definition in \eqref{eq20170319_4}, the matrix
$\mathbf{A}(\x)=(a_{kj})_{1\le j,k\le M}$ is symmetric positive
definite (\cite{xiu2009efficient}).
Thus,  equations \eqref{eq9}-\eqref{eq20160316_4} can be rewritten as the following:
\begin{align}\label{eq14}
\frac{\partial^2\widehat{\mathbf{v}}}{\partial t^2}(t,\x) &= \mathrm{div}(\mathbf{A}(\x)\widehat{\mathbf{S}}(t,\x)),\\
\widehat{\mathbf{S}}(t,\x) &= \mathbf{A}(\x)\nabla\widehat{\mathbf{v}}(t,\x),\label{eq15}
\end{align}
with initial and the boundary conditions
\begin{align}
\widehat{\mathbf{v}}(0,\x)&=\widehat{\mathbf{v}}_0(\x),\quad \widehat{\mathbf{v}}_t(0,\x)=\widehat{\mathbf{v}}_{1}(\x),\label{eq_add4}\\
\widehat{\mathbf{v}}(t,\x)|_{\partial\mathcal{D}}&=0.\label{eq_add1}
\end{align}



\subsection{LDG discretization}\label{sec-3}

To look for numerical approximation of \eqref{eq14}-\eqref{eq_add1},
we discretize the domain $\mathcal{D}$ into $K_{ij} := I_i\times J_j := [x_{i-\frac12},x_{i+\frac12}]\times[z_{j-\frac12},z_{j+\frac12}]$ for $1\leq i \leq N_x, 1\leq j\leq N_z$ and consider the following piecewise polynomial space 
\begin{align}
V_h^k := \bigl\{r\in L^2(\mathcal{D}): r|_{D_{ij}}\in P^k(K_{ij}), i =1,2, \cdots, N_x, j =1,2, \cdots, N_z\bigr\},\label{eq20170316_7}
\end{align}
where $P^k(K_{ij})$ denotes the space of polynomials with degree up to $k$ in the domain $K_{ij}$.
We define $\mathbf{V}_h^k$ as a space of vectored functions whose entries are in $V_h^k$.
In the following we use dot ($\cdot$) to denote a binary operation between two vectors or matrices which calculates the inner product of the corresponding row vectors (scalar multiplication in the case of vectors) and outputs a single column vector. The divergence operator is applied in a row-wise fashion.

The LDG method for Eqs.~(\ref{eq14})-(\ref{eq15}) is to seek $\widehat{\mathbf{v}}_h\in \mathbf{H}^2([0,T];\mathbf{V}_h^k)$, $\widehat{\mathbf{S}}_h\in (\mathbf{L}^2([0,T];\mathbf{V}_h^k))^2$ such that
\begin{alignat}{2}\label{eq19}
\int_{K_{ij}} \frac{\partial^2\widehat{\mathbf{v}}_h}{\partial t^2}\cdot\mathbf{p}_hd\x +  \int_{K_{ij}}\mathbf{A}\widehat{\mathbf{S}}_h\cdot\nabla\mathbf{p}_hd\x- (\widehat{\mathbf{A}\widehat{\mathbf{S}}}_h\cdot\bm{\nu},\mathbf{p}_h)_{\partial K_{ij}}& = 0\qquad \forall \mathbf{p}_h\in \mathbf{V}_h^k,\\
\int_{K_{ij}}\widehat{\mathbf{S}}_h\cdot\mathbf{w}_hd\x +  \int_{K_{ij}}\mathbf{A}\widehat{\mathbf{v}}_h\cdot\div(\mathbf{w}_h)d\x +  \int_{K_{ij}}\bar{\mathbf{A}}\widehat{\mathbf{v}}_h\cdot\mathbf{w}_hd\x&\label{eq20}\\
\qquad-(\widehat{\mathbf{A}\widehat{\mathbf{v}}_h},\mathbf{w}_h\cdot\bm{\nu})_{\partial K_{ij}} & = 0\qquad \forall \mathbf{w}_h\in (\mathbf{V}_h^k)^2,\notag
\end{alignat}
subject to the initial conditions $\widehat{\mathbf{v}}_h(0,\x) = P^+_h \widehat{\mathbf{v}}_0(\x), (\widehat{\mathbf{v}}_h)_t(0,\x) = P_h \widehat{\mathbf{v}}_{1}(\x)$, where the projections $P^+_h$ and $P_h$ will be specified later in Section \ref{sec-4}. In Eq.~(\ref{eq20}), $\bar{\mathbf{A}}$ denotes the matrix with each entry being the gradient of the corresponding entry of $\mathbf{A}$.

A critical step is to choose the numerical fluxes, which ultimately determines the property of the resulting scheme.
Assuming that $\mathbf{A}$ is piecewise smooth and the possible discontinuity occurs only along the direction aligned with the spatial discretization. We choose the flux associated with $\mathbf{A}$ to be the same as the test functions, namely, from inside of the cell in \eqref{eq20},  then \eqref{eq20} becomes
\begin{align}\label{eq21}
\int_{K_{ij}}\widehat{\mathbf{S}}_h\cdot\mathbf{w}_hd\x +  \int_{K_{ij}}\mathbf{A}\widehat{\mathbf{v}}_h\cdot\div(\mathbf{w}_h)d\x +  \int_{K_{ij}}\bar{\mathbf{A}}\widehat{\mathbf{v}}_h\cdot\mathbf{w}_hd\x&\\
\qquad-(\mathbf{A}\widehat{\widehat{\mathbf{v}}_h},\mathbf{w}_h\cdot\bm{\nu})_{\partial K_{ij}} & = 0\qquad \forall \mathbf{w}_h\in (\mathbf{V}_h^k)^2.\notag
\end{align}


Writing more explicitly, the LDG method \eqref{eq19} and \eqref{eq21} is to seek $\widehat{\mathbf{v}}_h\in \mathbf{H}^2([0,T];\mathbf{V}_h^k)$, $\widehat{\mathbf{S}}_h\in (\mathbf{L}^2([0,T];\mathbf{V}_h^k))^2$ such that
\begin{alignat}{2}\label{eq20170317_1}
\int_{K_{ij}} \frac{\partial^2\widehat{\mathbf{v}}_h}{\partial t^2}\cdot\mathbf{p}_hd\x +  \int_{K_{ij}}\mathbf{A}\widehat{\mathbf{S}}_h\cdot\nabla\mathbf{p}_hd\x-(\widehat{\mathbf{A}\widehat{\mathbf{S}}^1_h},\mathbf{p}_h^-)_{J_{j}}+(\widehat{\mathbf{A}\widehat{\mathbf{S}}^1_h},\mathbf{p}_h^+)_{J_{j}}&\\
-(\widetilde{\mathbf{A}\widehat{\mathbf{S}}^2_h},\mathbf{p}_h^-)_{I_{i}}+(\widetilde{\mathbf{A}\widehat{\mathbf{S}}^2_h},\mathbf{p}_h^+)_{I_{i}}& = 0\qquad \forall \mathbf{p}_h\in \mathbf{V}_h^k,\notag\\
\int_{K_{ij}}\widehat{\mathbf{S}}^1_h\cdot\mathbf{w}_h^1d\x +  \int_{K_{ij}}\mathbf{A}\widehat{\mathbf{v}}_h\cdot(\mathbf{w}_h^1)_xd\x +  \int_{K_{ij}}\mathbf{A}_x\widehat{\mathbf{v}}_h\cdot\mathbf{w}_h^1d\x&\label{eq20170317_2}\\
\qquad-(\mathbf{A}\widehat{\widehat{\mathbf{v}}_h},(\mathbf{w}_h^1)^-)_{J_{j}}+(\mathbf{A}\widehat{\widehat{\mathbf{v}}_h},(\mathbf{w}_h^1)^+)_{J_{j}} & = 0\qquad \forall \mathbf{w}_h^1\in \mathbf{V}_h^k,\notag\\
\int_{K_{ij}}\widehat{\mathbf{S}}^2_h\cdot\mathbf{w}_h^2d\x +  \int_{K_{ij}}\mathbf{A}\widehat{\mathbf{v}}_h\cdot(\mathbf{w}_h^2)_yd\x +  \int_{K_{ij}}\mathbf{A}_y\widehat{\mathbf{v}}_h\cdot\mathbf{w}_h^2d\x&\label{eq20170317_3}\\
\qquad-(\mathbf{A}\widetilde{\widehat{\mathbf{v}}_h},(\mathbf{w}_h^2)^-)_{I_{i}}+(\mathbf{A}\widetilde{\widehat{\mathbf{v}}_h},(\mathbf{w}_h^2)^+)_{I_{i}} & = 0\qquad \forall \mathbf{w}_h^2\in \mathbf{V}_h^k,\notag
\end{alignat}
subject to the initial conditions $\widehat{\mathbf{v}}_h(0,\x) = P^+_h \widehat{\mathbf{v}}_0(\x), (\widehat{\mathbf{v}}_h)_t(0,\x) = P_h \widehat{\mathbf{v}}_{1}(\x)$. 
Here $\widehat{\mathbf{S}}^i_h$ denotes the $i$-th column of $\widehat{\mathbf{S}}_h$. In the boundary terms of \eqref{eq20170317_2}-\eqref{eq20170317_3} the matrix $\mathbf{A}$ will be evaluated from the inside of the cell as in \eqref{eq21}. As for the numerical fluxes in Eqs.~\eqref{eq20170317_1}--\eqref{eq20170317_3}, we choose alternating flux, that is,
\begin{align}\label{eq22}
\widehat{\mathbf{A}\widehat{\mathbf{S}}_h^1} = \mathbf{A}^-(\widehat{\mathbf{S}}_h^1)^-,\qquad \widehat{\widehat{\mathbf{v}}_h} = \widehat{\mathbf{v}}^+_h,
\end{align}
or
\begin{align}\label{eq23}
\widehat{\mathbf{A}\widehat{\mathbf{S}}_h^1} = \mathbf{A}^+(\widehat{\mathbf{S}}_h^1)^+,\qquad \widehat{\widehat{\mathbf{v}}_h} = \widehat{\mathbf{v}}^-_h,
\end{align}
where $\mathbf{A}^+$ and $\mathbf{A}^-$ denote the matrices obtained by choosing $a_{kj}^+$ and $a_{kj}^-$ as their $kj$-th compotents respectively for each $kj$-th component $a_{kj}$ of matrix $\mathbf{A}$. Similarly, we can choose
\begin{align}\label{eq20170317_5}
\widetilde{\mathbf{A}\widehat{\mathbf{S}}_h^2} = \mathbf{A}^-(\widehat{\mathbf{S}}_h^2)^-,\qquad \widetilde{\widehat{\mathbf{v}}_h} = \widehat{\mathbf{v}}^+_h,
\end{align}
or
\begin{align}\label{eq20170317_6}
\widetilde{\mathbf{A}\widehat{\mathbf{S}}_h^2} = \mathbf{A}^+(\widehat{\mathbf{S}}_h^2)^+,\qquad \widetilde{\widehat{\mathbf{v}}_h} = \widehat{\mathbf{v}}^-_h.
\end{align}

\subsection{Semi-discrete energy law}

Using the fluxes defined above, we can prove that the semi-discrete method in \eqref{eq19} and \eqref{eq21} is energy conserving.
Here we only consider the case in \eqref{eq22} and \eqref{eq20170317_5}, and the proof with \eqref{eq23} and \eqref{eq20170317_6} is similar. 

\begin{theorem}\label{thm1}
The semi-discretized energy
\begin{align}\label{eq24}
E_h(t) :=\int_{\mathcal{D}}
  \left(\frac{\partial\widehat{\mathbf{v}}_h}{\partial t} \cdot
  \frac{\partial\widehat{\mathbf{v}}_h}{\partial t}
  +\widehat{\mathbf{S}}_h\cdot \widehat{\mathbf{S}}_h \right) d\x
\end{align}
 is conserved by the semi-discretized scheme \eqref{eq19} and \eqref{eq21} for all time $t>0$.
\end{theorem}
\begin{proof}
By taking the time derivative of Eq.~\eqref{eq21} and choosing $\mathbf{w}_h = \widehat{\mathbf{S}}_h$, we obtain
\begin{align}\label{eq25}
\int_{K_{ij}}(\widehat{\mathbf{S}}_h)_t\cdot\widehat{\mathbf{S}}_hd\x +  \int_{K_{ij}}\mathbf{A}(\widehat{\mathbf{v}}_h)_t\cdot\div(\widehat{\mathbf{S}}_h)d\x +  \int_{K_{ij}}\bar{\mathbf{A}}(\widehat{\mathbf{v}}_h)_t\cdot\widehat{\mathbf{S}}_hd\x-(\mathbf{A}(\widehat{\mathbf{v}}_h^+)_t,\widehat{\mathbf{S}}_h\cdot\bm{\nu})_{\partial K_{ij}} = 0.
\end{align}
Taking $\mathbf{p}_h = (\widehat{\mathbf{v}}_h)_t$ in \eqref{eq19} yields
\begin{align}\label{eq26}
\int_{K_{ij}} \frac{\partial^2\widehat{\mathbf{v}}_h}{\partial t^2}\cdot(\widehat{\mathbf{v}}_h)_td\x +  \int_{K_{ij}}\mathbf{A}\widehat{\mathbf{S}}_h\cdot\nabla(\widehat{\mathbf{v}}_h)_td\x- (\mathbf{A}^-\widehat{\mathbf{S}}_h^-\cdot\bm{\nu},(\widehat{\mathbf{v}}_h)_t)_{\partial K_{ij}} = 0.
\end{align}
Adding \eqref{eq25} to \eqref{eq26} and using integration by parts on the second term of \eqref{eq26}, we have
\begin{align}\label{eq27}
\int_{K_{ij}}(\widehat{\mathbf{S}}_h)_t\cdot\widehat{\mathbf{S}}_hd\x+\int_{K_{ij}} \frac{\partial^2\widehat{\mathbf{v}}_h}{\partial t^2}\cdot(\widehat{\mathbf{v}}_h)_td\x+(\mathbf{A}\widehat{\mathbf{S}}_h\cdot\bm{\nu},(\widehat{\mathbf{v}}_h)_t)_{\partial K_{ij}}&\\
-(\mathbf{A}(\widehat{\mathbf{v}}_h^+)_t,\widehat{\mathbf{S}}_h\cdot\bm{\nu})_{\partial K_{ij}}- (\mathbf{A}^-\widehat{\mathbf{S}}_h^-\cdot\bm{\nu},(\widehat{\mathbf{v}}_h)_t)_{\partial K_{ij}}&=0.\notag
\end{align}
After summing over $K_{ij}$, Eq.~\eqref{eq27} can be written as
\begin{align}\label{eq28}
&\sum_{K_{ij}\in\mathcal{T}_h}\int_{K_{ij}}(\widehat{\mathbf{S}}_h)_t\cdot\widehat{\mathbf{S}}_hd\x+\sum_{K_{ij}\in\mathcal{T}_h}\int_{K_{ij}} \frac{\partial^2\widehat{\mathbf{v}}_h}{\partial t^2}\cdot(\widehat{\mathbf{v}}_h)_td\x+\sum_{E\in\mathcal{E}_h}(\mathbf{A}^-\widehat{\mathbf{S}}_h^-\cdot\bm{\nu},(\widehat{\mathbf{v}}^-_h)_t)_{E}\\
&\qquad-\sum_{E\in\mathcal{E}_h}(\mathbf{A}^+\widehat{\mathbf{S}}_h^+\cdot\bm{\nu},(\widehat{\mathbf{v}}^+_h)_t)_{E}-\sum_{E\in\mathcal{E}_h}(\mathbf{A}^-(\widehat{\mathbf{v}}_h^+)_t,\widehat{\mathbf{S}}^-_h\cdot\bm{\nu})_{E}+\sum_{E\in\mathcal{E}_h}(\mathbf{A}^+(\widehat{\mathbf{v}}_h^+)_t,\widehat{\mathbf{S}}^+_h\cdot\bm{\nu})_{E}\notag\\
&\qquad-\sum_{E\in\mathcal{E}_h}(\mathbf{A}^-\widehat{\mathbf{S}}_h^-\cdot\bm{\nu},(\widehat{\mathbf{v}}_h^-)_t )_{E}+\sum_{E\in\mathcal{E}_h}(\mathbf{A}^-\widehat{\mathbf{S}}_h^-\cdot\bm{\nu},(\widehat{\mathbf{v}}_h^+)_t )_{E}=0.\notag
\end{align}
By applying Dirichlet boundary conditions \eqref{eq_add1} and summing over $K_{ij}$, we get
\begin{align}\label{eq29}
\frac{d}{dt}\int_{\mathcal{D}}
  \left(\frac{\partial\widehat{\mathbf{v}}_h}{\partial t} \cdot
  \frac{\partial\widehat{\mathbf{v}}_h}{\partial t}
  +\widehat{\mathbf{S}}_h\cdot \widehat{\mathbf{S}}_h \right) d\x=0.
\end{align}

Therefore, $E_h(t)$ is invariant in time.
\end{proof}

\subsection{Fully discrete energy law}
Next, we consider the fully-discrete LDG method with leap-frog time discretization. Let $0= t_0\le t_1\le\cdots\le t_N=T$ be a uniform partition of the interval $[0,T]$ with time step size $\Delta t$. We use $\widehat{\mathbf{v}}_h^{n}, \widehat{\mathbf{S}}_h^{n}$ to denote the numerical solutions at $t=t_n$.
Thus the scheme is to seek $\widehat{\mathbf{v}}_h^{n+1}\in \mathbf{V}_h^k$, $\widehat{\mathbf{S}}_h^n\in (\mathbf{V}_h^k)^2$ such that for all $K_{ij}$, the following equations hold:
\begin{align}\label{eq30}
\int_{K_{ij}} \frac{\widehat{\mathbf{v}}_h^{n+1}-2\widehat{\mathbf{v}}_h^{n}+\widehat{\mathbf{v}}_h^{n-1}}{(\Delta t)^2}\cdot\mathbf{p}_hd\x + \int_{K_{ij}}\mathbf{A}\widehat{\mathbf{S}}^n_h\cdot\nabla\mathbf{p}_hd\x- (\mathbf{A}^-(\widehat{\mathbf{S}}_h^n)^-\cdot\bm{\nu},\mathbf{p}_h)_{\partial K_{ij}} & = 0\qquad \forall \mathbf{p}_h\in \mathbf{V}_h^k,\\
\int_{K_{ij}}\widehat{\mathbf{S}}^n_h\cdot\mathbf{w}_hd\x +  \int_{K_{ij}}\mathbf{A}\widehat{\mathbf{v}}_h^{n}\cdot\div(\mathbf{w}_h)d\x +  \int_{K_{ij}}\bar{\mathbf{A}}\widehat{\mathbf{v}}^{n}_h\cdot\mathbf{w}_hd\x&\label{eq31}\\
\qquad-(\mathbf{A}(\widehat{\mathbf{v}}^n_h)^+,\mathbf{w}_h\cdot\bm{\nu})_{\partial K_{ij}} & = 0\qquad \forall \mathbf{w}_h\in (\mathbf{V}_h^k)^2,\notag
\end{align}
subject to the initial conditions $\widehat{\mathbf{v}}_h^0(0,\x) = P^+_h \widehat{\mathbf{v}}_0(\x), (\widehat{\mathbf{v}}_h)_t^0(0,\x) = P_h \widehat{\mathbf{v}}_{00}(\x)$. 

%

In the following we show the fully-discrete energy law.
\begin{theorem}\label{thm2}
The fully-discrete energy, defined by
\begin{align}\label{eq35}
E_h^{n+1} :=\bigg\|\frac{\widehat{\mathbf{v}}_h^{n+1}-\widehat{\mathbf{v}}_h^{n}}{\Delta t}\bigg\|^2+\bigg\|\frac{\widehat{\mathbf{S}}_h^{n+1}+\widehat{\mathbf{S}}_h^{n}}{2}\bigg\|^2-\frac{(\Delta t)^2)}{4}\bigg\|\frac{\widehat{\mathbf{S}}_h^{n+1}-\widehat{\mathbf{S}}_h^{n}}{\Delta t}\bigg\|^2
\end{align}
 is conserved by the fully-discrete scheme \eqref{eq30} and \eqref{eq31} for all $n$.
\end{theorem}
\begin{proof}
In \eqref{eq30}, we choose the test function to be $\mathbf{p}_h = \frac{\widehat{\mathbf{v}}_h^{n+1}-\widehat{\mathbf{v}}_h^{n-1}}{2\Delta t}$, then
\begin{align}\label{eq36}
\int_{K_{ij}} \frac{\widehat{\mathbf{v}}_h^{n+1}-2\widehat{\mathbf{v}}_h^{n}+\widehat{\mathbf{v}}_h^{n-1}}{(\Delta t)^2}\cdot\frac{\widehat{\mathbf{v}}_h^{n+1}-\widehat{\mathbf{v}}_h^{n-1}}{2\Delta t}d\x +  \int_{K_{ij}}\mathbf{A}\widehat{\mathbf{S}}^n_h\cdot\nabla(\frac{\widehat{\mathbf{v}}_h^{n+1}-\widehat{\mathbf{v}}_h^{n-1}}{2\Delta t})d\x&\\
- (\mathbf{A}^-(\widehat{\mathbf{S}}_h^n)^-\cdot\bm{\nu},\frac{\widehat{\mathbf{v}}_h^{n+1}-\widehat{\mathbf{v}}_h^{n-1}}{2\Delta t})_{\partial K_{ij}} &= 0.\notag
\end{align}

Considering the equation \eqref{eq31} at time $t_{n-1}$ and $t_{n+1}$, and taking the test function $\mathbf{w}_h = \frac{1}{2\Delta t}\widehat{\mathbf{S}}_h^n$, we obtain
\begin{align}\label{eq37}
\int_{K_{ij}}\frac{\widehat{\mathbf{S}}_h^{n+1}-\widehat{\mathbf{S}}_h^{n-1}}{2\Delta t}\cdot\widehat{\mathbf{S}}_h^nd\x +  \int_{K_{ij}}\mathbf{A}\frac{\widehat{\mathbf{v}}_h^{n+1}-\widehat{\mathbf{v}}_h^{n-1}}{2\Delta t}\cdot\div(\widehat{\mathbf{S}}_h^n)d\x+\int_{K_{ij}}\bar{\mathbf{A}}\frac{\widehat{\mathbf{v}}_h^{n+1}-\widehat{\mathbf{v}}_h^{n-1}}{2\Delta t}\cdot\widehat{\mathbf{S}}_h^nd\x&\\
-(\mathbf{A}\bigl(\frac{\widehat{\mathbf{v}}_h^{n+1}-\widehat{\mathbf{v}}_h^{n-1}}{2\Delta t}\bigr)^+,\widehat{\mathbf{S}}_h^n\cdot\bm{\nu})_{\partial K_{ij}} &= 0.\notag
\end{align}
By adding \eqref{eq36} to \eqref{eq37}, summing over $K_{ij}$, and using integration by parts, we have
\begin{align}\label{eq38}
0&=\sum_{K_{ij}\in\mathcal{T}_h}\int_{K_{ij}} \frac{\widehat{\mathbf{v}}_h^{n+1}-2\widehat{\mathbf{v}}_h^{n}+\widehat{\mathbf{v}}_h^{n-1}}{(\Delta t)^2}\cdot\frac{\widehat{\mathbf{v}}_h^{n+1}-\widehat{\mathbf{v}}_h^{n-1}}{2\Delta t}d\x+\sum_{K_{ij}\in\mathcal{T}_h}\int_{K_{ij}}\frac{\widehat{\mathbf{S}}_h^{n+1}-\widehat{\mathbf{S}}_h^{n-1}}{2\Delta t}\cdot\widehat{\mathbf{S}}_h^nd\x\\
&=\sum_{K_{ij}\in\mathcal{T}_h}\int_{K_{ij}}\frac{(\widehat{\mathbf{v}}_h^{n+1}-\widehat{\mathbf{v}}_h^{n})-(\widehat{\mathbf{v}}_h^{n}-\widehat{\mathbf{v}}_h^{n-1})}{(\Delta t)^2}\cdot\frac{(\widehat{\mathbf{v}}_h^{n+1}-\widehat{\mathbf{v}}_h^{n})+(\widehat{\mathbf{v}}_h^{n}-\widehat{\mathbf{v}}_h^{n-1})}{2\Delta t}d\x\notag\\
&+\sum_{K_{ij}\in\mathcal{T}_h}\int_{K_{ij}}\frac{\widehat{\mathbf{S}}_h^{n+1}+2\widehat{\mathbf{S}}_h^{n}+\widehat{\mathbf{S}}_h^{n-1}}{4}\cdot\frac{\widehat{\mathbf{S}}_h^{n+1}-\widehat{\mathbf{S}}_h^{n-1}}{2\Delta t}d\x-\sum_{K_{ij}\in\mathcal{T}_h}\int_{K_{ij}}\cdot\frac{\widehat{\mathbf{S}}_h^{n+1}-2\widehat{\mathbf{S}}_h^{n}+\widehat{\mathbf{S}}_h^{n-1}}{4}\frac{\widehat{\mathbf{S}}_h^{n+1}-\widehat{\mathbf{S}}_h^{n-1}}{2\Delta t}d\x\notag\\
&=\  \frac{1}{2\Delta t}(E_h^{n+1}-E_h^{n}),\notag
\end{align}
with $E_h^{n}$ defined in \eqref{eq35}. Thus the discrete energy is conserved over time.
\end{proof}

\begin{remark}
There is a term with uncertain sign in $E_h^{n+1}$, and this term comes from the use of the explicit leapfrog scheme. By some calculations, we know
\begin{align*}
E_h^{n+1} =\bigg\|\frac{\widehat{\mathbf{v}}_h^{n+1}-\widehat{\mathbf{v}}_h^{n}}{\Delta t}\bigg\|^2+(\widehat{\mathbf{S}}_h^n,\widehat{\mathbf{S}}_h^{n+1}).
\end{align*}
Formally $\Delta t$ needs to be small enough to guarantee $E_h^{n+1}\ge0$.
\end{remark}

\section{Error estimates}\label{sec-4}
In this section, we provide error estimate for the spatial discretization in the semi-discrete scheme \eqref{eq19} and \eqref{eq21}.
We will show that the error bound is optimal and is linear in time. 
Let $u(t,\x,\y)$ and $\q(t,\x,\y)$ be the exact solution of \eqref{eq_u} and \eqref{eq_q}, and $u_h(t,\x,\y)$ and $\q_h(t,\x,\y)$ are numerical solutions
\begin{align}
u_h(t,\x,\y) = \sum_{m=1}^M(\widehat{v}_m)_h(t,\x)\Phi_m(\y),\notag\\
\q_h(t,\x,\y) = \sum_{m=1}^M(\widehat{\r}_m)_h(t,\x)\Phi_m(\y),\notag\\
\end{align}
where $(\widehat{v}_m)_h$ and $(\widehat{\r}_m)_h$ are the $m$-th row of $\widehat{\mathbf{v}}_h$ and $\widehat{\mathbf{S}}_h$.
We consider the errors:
\begin{eqnarray}\label{eq20161222_add1}
e_u &=& u - u_h = (u - u_M) + (u_M - u_h) \label{eq20161222_add1}\\
e_\q &=& \q - \q_h = (\q - \q_M) + (\q_M - \q_h) ,\label{error_q}
\end{eqnarray}
where the $u_M$ and $\q_M$ are the gPC approximations defined in \eqref{eq5} and \eqref{eq20160316_1}. 
We call the first term on the right-hand side of \eqref{eq20161222_add1}-\eqref{error_q} the gPC approximation error and the second term the spatial discretization error. 
In the following, we provide the error estimates in semi-discrete energy norm and show that the convergence is optimal.
\begin{theorem}\label{thm20161222_add1}
Let $e_u$ and $e_\q$ defined by \eqref{eq20161222_add1} and \eqref{error_q}, and initial conditions satisfy 
\begin{align}\label{eq20180131}
\widehat{\mathbf{v}}_h(\x,0)=P_h^+\widehat{\mathbf{v}}(\x,0),\qquad(\widehat{\mathbf{v}}_h)_t(\x,0)=P_h\widehat{\mathbf{v}}_t(\x,0).
\end{align}
For any given $\epsilon_M$, if we choose $M$ in \eqref{eq5}-\eqref{eq20160316_1} sufficiently large so that 
\begin{equation}\label{eq20170319_11}
\sum_{j=M+1}^{\infty}\|\widehat{\mathbf{p}}_j\|_{H^1(\D)}\le \epsilon_M,\quad
\sum_{j=M+1}^{\infty}\|\widehat{v}_j\|_{H^1(\D)}\le \epsilon_M,
\end{equation}
then with \eqref{eq5}-\eqref{eq20160316_1} and the LDG approximation \eqref{eq19} and \eqref{eq21}, the error estimate in the energy norm is 
\begin{equation}
\bigl(\mathbb{E}[\|(e_u)_t\|^2_{L^2(\mathcal{D})}]\bigr)^{1\slash2}+\bigl(\mathbb{E}[\|e_\q\|^2_{L^2(\mathcal{D})}]\bigr)^{1\slash2}
\le C(t+1)\epsilon_M+ C(t+1)h^{k+1}.\notag
\end{equation}
\end{theorem}

\begin{proof}
We divide our proof into two parts, corresponding to bounds for the gPC approximation error and semi-discretization error, respectively.\\

\noindent{\bf{Part 1 (The gPC approximation error).}} First we rewrite Eqs.~\eqref{eq20170319_2}-\eqref{eq20170319_3} as
\begin{align}\label{eq60}
\frac{\partial^2 v_k}{\partial t^2}(t,x)&=\sum_{j=1}^{M}\mathrm{div}(a_{kj}(\x) \r_j)+\sum_{j=M+1}^{\infty}\mathrm{div}(a_{kj}(\x) \r_j),\quad k=1,2,\cdots,\\
\r_k(t,x)&=\sum_{j=1}^{M}a_{kj}(\x) \nabla v_j+\sum_{j=M+1}^{\infty}a_{kj}(\x) \nabla v_j.\label{eq20170319_6}
\end{align}
Denoting $\mathbf{v} = (v_1,v_2,\ldots,v_M)^T$ and $\mathbf{S} = (\p_1^T,\p_2^T,\ldots,\p_M^T)^T$, then Eqs~\eqref{eq60}-\eqref{eq20170319_6} for $k=1,\cdots, M$ can be written as
\begin{align}
\frac{\partial^2 \mathbf{v}(t,\x)}{\partial t^2} &=  \div(\mathbf{A}(\x)\mathbf{S}(t,\x))+\mathbf{r}(t,\x),\label{eq58}\\
\mathbf{S}(t,\x) &=  \mathbf{A}(\x)\nabla\mathbf{v}(t,\x)+{\mathbf{R}}(t,\x),\label{eq20170319_7}
\end{align}
where $\mathbf{A}$ is defined as in \eqref{eq20170319_4}.
In \eqref{eq58}, $\mathbf{r}(t,\x)$ is a vector, with the $k$-th component defined by
\[
\mathbf{r}_k = \sum_{j=M+1}^{\infty}\div(a_{jk}(\x)\r_j),
\]
and in \eqref{eq20170319_7}, ${\mathbf{R}}(t,\x)$ is a matrix with its $k$-th row as
\[
\mathbf{R}_k = \sum_{j=M+1}^{\infty}a_{jk}(\x)\nabla v_j.
\]
Subtracting Eqs.~\eqref{eq14}-\eqref{eq15} from Eqs.~\eqref{eq58}--\eqref{eq20170319_7}, we get
\begin{align}
\frac{\partial^2 (\mathbf{v}-\mathbf{\widehat{v}})}{\partial t^2} &=  \div(\mathbf{A}(\x)(\mathbf{S}-\mathbf{\widehat{S}}))+\mathbf{r},\label{eq64}\\
\mathbf{S}-\widehat{\mathbf{S}}&=\mathbf{A}(\x)\nabla(\mathbf{v}-\mathbf{\widehat{v}})+{\mathbf{R}}.\label{eq20170319_10}
\end{align}

We first multiply \eqref{eq64} by $\mathbf{v}_t-\widehat{\mathbf{v}}_t$ and integrate in space over $\mathcal{D}$, and then take the time derivative of \eqref{eq20170319_10}, followed by multiplying \eqref{eq20170319_10} with $\mathbf{S}_t-\widehat{\mathbf{S}}_t$ and integration over $\mathcal{D}$. With the fact that the coefficients $a_{jk}$ are bounded, we obtain the estimate:
\begin{align}\
&\qquad\frac12\frac{\partial}{\partial t}\big(\| \mathbf{v}_t-\widehat{\mathbf{v}}_t\|_{L^2(\mathcal{D})}^2+\|(\mathbf{S}-\widehat{\mathbf{S}})\|_{L^2(\mathcal{D})}^2\bigr)\notag\\
&=(\mathbf{r},\mathbf{v}_t-\widehat{\mathbf{v}}_t)+({\mathbf{R}},\mathbf{S}-\widehat{\mathbf{S}})\notag\\
&\le\|\mathbf{r}\|_{L^2(\mathcal{D})} \|\mathbf{v}_t-\widehat{\mathbf{v}}_t\|_{L^2(\mathcal{D})}+\|{\mathbf{R}}\|_{L^2(\mathcal{D})} \|\mathbf{S}-\widehat{\mathbf{S}}\|_{L^2(\mathcal{D})}\notag\\
&\le C\sum_{j=M+1}^{\infty}\|\mathbf{p}_j\|_{H^1(\D)}\|\mathbf{v}_t-\widehat{\mathbf{v}}_t\|_{L^2(\mathcal{D})}+C\sum_{j=M+1}^{\infty}\|v_j\|_{H^1(\mathcal{D})}\|\mathbf{S}-\widehat{\mathbf{S}}\|_{L^2(\mathcal{D})}.\label{eq63}
\end{align}
By \eqref{eq20170319_11}, we have
\begin{align}\label{eq20181024_1}
\left(\mathbb{E}\left[\|(u)_{t}-(u_M)_{t}\|_{L^2(\mathcal{D})}^2\right]\right)^{1\slash2} + \left(\mathbb{E}\left[\|\q-\q_M\|^2_{L^2(\mathcal{D})}\right]\right)^{1\slash2} \le C(t+1)\epsilon_M.
\end{align}

\noindent{\bf{Part 2 (The spatial discretization error).}} 
Consider the weak formulation of \eqref{eq14}-\eqref{eq15}: finding $\widehat{\mathbf{v}}\in \mathbf{H}^2(\mathcal{T};\mathbf{H}^{-1}( \mathcal{D}))\cap \mathbf{L}^2(\mathcal{T};\mathbf{H}^1(\mathcal{D})),$ $\widehat{\mathbf{S}}\in (\mathbf{L}^2(\mathcal{T};\mathbf{H}^1( \mathcal{D})))^2$ such that
\begin{align}\label{eq20161218_3}
\int_{K_{ij}} \frac{\partial^2\widehat{\mathbf{v}}}{\partial t^2}\cdot\mathbf{p}d\x +  \int_{K_{ij}}\mathbf{A}\widehat{\mathbf{S}}\cdot\nabla\mathbf{p}d\x- (\mathbf{A}^-(\widehat{\mathbf{S}})^-\cdot\bm{\nu},\mathbf{p})_{\partial{K_{ij}}} & = 0\qquad \forall \mathbf{p}\in \mathbf{H}^1( \mathcal{D}),\\
\int_{K_{ij}}\widehat{\mathbf{S}}\cdot\mathbf{w}d\x +  \int_{K_{ij}}\mathbf{A}\widehat{\mathbf{v}}\cdot\div\mathbf{w}d\x +  \int_{K_{ij}}\bar{\mathbf{A}}\widehat{\mathbf{v}}\cdot\mathbf{w}d\x\label{eq20161218_4}\\
-(\mathbf{A}\widehat{\mathbf{v}}^+,\mathbf{w}\cdot\bm{\nu})_{\partial{K_{ij}}} &= 0\qquad \forall \mathbf{w}\in (\mathbf{H}^1(\mathcal{D}))^2.\notag
\end{align}
Note that the jump conditions $\widehat{\mathbf{v}}^-=\widehat{\mathbf{v}}^+$ and $A^+(\widehat{\mathbf{S}})^+=A^-(\widehat{\mathbf{S}})^-$ are assumed on the mesh boundaries, and $\bar{\mathbf{A}}$ is defined in Section \ref{sec-3}.

On the other hand, the LDG approximation is to look for $\widehat{\mathbf{v}}_h$ and $\widehat{\mathbf{S}}_h$ such that 
\begin{align}\label{eq20161218_5}
\int_{K_{ij}} \frac{\partial^2\widehat{\mathbf{v}}_h}{\partial t^2}\cdot\mathbf{p}_hd\x +  \int_{K_{ij}}\mathbf{A}\widehat{\mathbf{S}}_h\cdot\nabla\mathbf{p}_hd\x- (\mathbf{A}^-(\widehat{\mathbf{S}}_h)^-\cdot\bm{\nu},\mathbf{p}_h)_{\partial{K_{ij}}} & = 0\qquad \forall \mathbf{p}_h\in \mathbf{V}_h^k,\\
\int_{K_{ij}}\widehat{\mathbf{S}}_h\cdot\mathbf{w}_hd\x +  \int_{K_{ij}}\mathbf{A}\widehat{\mathbf{v}}_h\cdot\div\mathbf{w}_hd\x +  \int_{K_{ij}}\bar{\mathbf{A}}\widehat{\mathbf{v}}_h\cdot\mathbf{w}_hd\x\label{eq20161218_4_add}\\
-(\mathbf{A}(\widehat{\mathbf{v}}_h)^+,\mathbf{w}_h\cdot\bm{\nu})_{\partial{K_{ij}}} &= 0\qquad \forall \mathbf{w}_h\in (\mathbf{V}_h^k)^2.\notag
\end{align}

Here we define $P_h$ to be the usual projection of a vectored function $\mathbf{u}$ associated with matrix $A$, that is,
\begin{align*}
(P_h\mathbf{u},\mathbf{Av})_{K_{ij}} = (\mathbf{u},\mathbf{Av})_{K_{ij}}\qquad\forall {\mathbf{v}}\in\mathbf{V}_h^k,
\end{align*}
and define $P_x^+$, $P_x^-$, $P_y^+$ and $P_y^-$ as the following special projections
\begin{alignat}{3}
(P_x^-\mathbf{u},\mathbf{Av})_{K_{ij}} &= (\mathbf{u},\mathbf{Av})_{K_{ij}},&&\quad\forall {\mathbf{v}}\in\mathbf{V}_h^{k-1}\quad\text{and}\quad(P_x^-\mathbf{u})^-(x_{i+\frac12}) = \mathbf{u}^-(x_{i+\frac12}),\notag\\
(P_x^+\mathbf{u},\mathbf{Av})_{K_{ij}} &= (\mathbf{u},\mathbf{Av})_{K_{ij}},&&\quad\forall {\mathbf{v}}\in\mathbf{V}_h^{k-1}\quad\text{and}\quad(P_x^+\mathbf{u})^+(x_{i-\frac12}) = \mathbf{u}^+(x_{i-\frac12}),\notag\\
(P_y^-\mathbf{u},\mathbf{Av})_{K_{ij}} &= (\mathbf{u},\mathbf{Av})_{K_{ij}},&&\quad\forall {\mathbf{v}}\in\mathbf{V}_h^{k-1}\quad\text{and}\quad(P_y^-\mathbf{u})^-(y_{i+\frac12}) = \mathbf{u}^-(y_{i+\frac12}),\notag\\
(P_y^+\mathbf{u},\mathbf{Av})_{K_{ij}} &= (\mathbf{u},\mathbf{Av})_{K_{ij}},&&\quad\forall {\mathbf{v}}\in\mathbf{V}_h^{k-1}\quad\text{and}\quad(P_y^+\mathbf{u})^+(y_{i-\frac12}) = \mathbf{u}^+(y_{i-\frac12}).\notag
\end{alignat}
We further define the errors by
\begin{alignat}{3}
\bar{e}_u &= \widehat{\mathbf{v}} - \widehat{\mathbf{v}}_h,\qquad \mathbf{\xi}_u &= \widehat{\mathbf{v}}-P_h^+\widehat{\mathbf{v}}, \qquad \eta_u&=P_h^+\widehat{\mathbf{v}}-\widehat{\mathbf{v}}_h,\notag\\
\bar{e}_{\q} &= \widehat{\mathbf{S}} - \widehat{\mathbf{S}}_h,\qquad \mathbf{\xi}_{\q} &= \widehat{\mathbf{S}}-P_h^-\widehat{\mathbf{S}}, \qquad \mathbf{\eta}_{\q}&=P_h^-\widehat{\mathbf{S}}-\widehat{\mathbf{S}}_h,\notag
\end{alignat}
where $P_h^+=P_x^+\otimes P_y^+$ and $P_h^-=P_x^-\otimes P_y^-$.\\

Subtracting \eqref{eq20161218_5}-\eqref{eq20161218_4_add} from \eqref{eq20161218_3}-\eqref{eq20161218_4},
and using the above definitions, we can rewrite the error equations into
\begin{align}\label{eq20161218_7}
\int_{K_{ij}} \frac{\partial^2\mathbf{\eta}_u}{\partial t^2}\cdot\mathbf{p}_hd\x +\int_{K_{ij}} \frac{\partial^2\mathbf{\xi}_u}{\partial t^2}\cdot\mathbf{p}_hd\x+ \int_{K_{ij}}\mathbf{\eta}_{\q}\cdot\mathbf{A}\nabla\mathbf{p}_hd\x&\\
- (\mathbf{\eta}_{\mathbf{q}}^-\cdot\bm{\nu},\mathbf{A}^-\mathbf{p}_h)_{\partial{K_{ij}}}&= 0\qquad \forall \mathbf{p}_h\in \mathbf{V}_h^k,\notag\\
\int_{K_{ij}}\mathbf{\xi}_{\q}\cdot\mathbf{w}_hd\x+\int_{K_{ij}}\mathbf{\eta}_{\q}\cdot\mathbf{w}_hd\x +  \int_{K_{ij}}\mathbf{\eta}_u\cdot\mathbf{A}\div\mathbf{w}_hd\x+  \int_{K_{ij}}\mathbf{\xi}_u\cdot\mathbf{A}\div\mathbf{w}_hd\x \label{eq20161218_8}\\
+  \int_{K_{ij}}\bar{\mathbf{A}}\mathbf{\eta}_u\cdot\mathbf{w}_hd\x +\int_{K_{ij}}\bar{\mathbf{A}}\mathbf{\xi}_u\cdot\mathbf{w}_hd\x\notag\\
-(\mathbf{A}\mathbf{\eta}_u^+,\mathbf{w}_h\cdot\bm{\nu})_{\partial{K_{ij}}}-(\mathbf{A}\mathbf{\xi}_u^+,\mathbf{w}_h\cdot\bm{\nu})_{\partial{K_{ij}}}&= 0\qquad \forall \mathbf{w}_h\in (\mathbf{V}_h^k)^2.\notag
\end{align}
Taking the time derivative of \eqref{eq20161218_8} and choosing $\mathbf{w}_h=\mathbf{\eta}_{\q}$ and $\mathbf{p}_h=(\mathbf{\eta}_u)_t$, 
the sum of these equations yields
\begin{align}\label{eq20161222_1}
&\qquad\int_{K_{ij}}(\mathbf{\eta}_u)_{tt}\cdot(\mathbf{\eta}_u)_{t}d\x+\int_{K_{ij}}(\mathbf{\eta}_{\q})_t\cdot\mathbf{\eta}_{\q}d\x\\
&=-\int_{K_{ij}}(\mathbf{\xi}_u)_{tt}\cdot(\mathbf{\eta}_u)_td\x- \int_{K_{ij}}\mathbf{\eta}_{\q}\cdot\mathbf{A}\nabla(\mathbf{\eta}_u)_{t}d\x-\int_{K_{ij}}(\mathbf{\xi}_{\q})_t\cdot\mathbf{\eta}_{\q}d\x \notag\\
&\quad-\int_{K_{ij}}(\mathbf{\eta}_u)_t\cdot\mathbf{A}\div\mathbf{\eta}_{\q}d\x-\int_{K_{ij}}(\mathbf{\xi}_u)_t\cdot\mathbf{A}\div\mathbf{\eta}_{\q}d\x-\int_{K_{ij}}\bar{\mathbf{A}}(\mathbf{\eta}_u)_t\cdot\mathbf{\eta}_{\q}d\x-\int_{K_{ij}}\bar{\mathbf{A}}(\mathbf{\xi}_u)_t\cdot\mathbf{\eta}_{\q}d\x\notag\\
&\quad + (\mathbf{\eta}_{\mathbf{q}}^-\cdot\bm{\nu},\mathbf{A}^-(\mathbf{\eta}_u)_{t})_{\partial{K_{ij}}}+(\mathbf{A}(\mathbf{\eta}_u^+)_t,\mathbf{\eta}_{\mathbf{q}}\cdot\bm{\nu})_{\partial{K_{ij}}}+(\mathbf{A}(\mathbf{\xi}_u^+)_t,\mathbf{\eta}_{\mathbf{q}}\cdot\bm{\nu})_{\partial{K_{ij}}}.\notag
\end{align}
By integration by parts to the fourth term on the right-hand side of \eqref{eq20161222_1}, and summing over all cells $K_{ij}$, we have
\begin{align}\label{eq20161222_2}
&\qquad\int_{\mathcal{T}_h}(\mathbf{\eta}_u)_{tt}\cdot(\mathbf{\eta}_u)_{t}d\x+\int_{\mathcal{T}_h}(\mathbf{\eta}_{\q})_t\cdot\mathbf{\eta}_{\q}d\x\\
&=-\int_{\mathcal{T}_h}(\mathbf{\xi}_u)_{tt}\cdot(\mathbf{\eta}_u)_td\x-\int_{\mathcal{T}_h}(\mathbf{\xi}_{\q})_t\cdot\mathbf{\eta}_{\q}d\x\notag\\
&\quad-\int_{K_{ij}}(\mathbf{\xi}_u)_t\cdot\mathbf{A}\div\mathbf{\eta}_{\q}d\x+(\mathbf{A}(\mathbf{\xi}_u^+)_t,\mathbf{\eta}_{\mathbf{q}}\cdot\bm{\nu})_{\partial{K_{ij}}}-\int_{K_{ij}}\bar{\mathbf{A}}(\mathbf{\xi}_u)_t\cdot\mathbf{\eta}_{\q}d\x.\notag
\end{align}
By the Cauchy-Schwarz's inequality and (3.3) in \cite{chou2014optimal} or Lemma 3.7 in \cite{dong2009analysis}, we have
\begin{align}\label{eq20161220_1}
&\qquad\frac12\frac{d}{dt}\bigl(\|(\mathbf{\eta}_u)_t\|^2+\|\mathbf{\eta}_{\q}\|^2\bigr)\\
&\le Ch^{k+1}\bigl(\|(\mathbf{\eta}_u)_t\|+\|\mathbf{\eta}_{\q}\|\bigr)+Ch^{k+1}\|u_t\|_{H^{k+2}}\|\mathbf{\eta}_{\q}\|+\|\bar{\mathbf{A}}\|_{L^{\infty}(\D)}\|(\xi_u)_t\|\|\mathbf{\eta}_{\q}\|\notag\\
&\le Ch^{k+1}\bigl(\|(\mathbf{\eta}_u)_t\|^2+\|\mathbf{\eta}_{\q}\|^2\bigr)^{1\slash2}.\notag
\end{align}
If we choose the initial conditions specifically to be \eqref{eq20180131}
then we have (\cite{chou2014optimal, xing2013energy})
\begin{align}\label{eq20171215}
\|(\mathbf{\eta}_u)_t(0)\|\le Ch^{k+1},\qquad\|\mathbf{\eta}_{\q}(0)\|\le Ch^{k+1},
\end{align}
and therefore
\begin{align}\label{eq20161220_2}
\bigl(\|(\mathbf{\eta}_u)_t\|^2+\|\mathbf{\eta}_{\q}\|^2\bigr)^{1\slash2}\le C(t+1)h^{k+1}.
\end{align}
By the properties of the projections, 
\begin{align}\label{eq20161222_10}
\bigl(\|(\bar{e}_u)_t\|^2+\|\bar{e}_{\q}\|^2\bigr)^{1\slash2}\le C(t+1)h^{k+1}.
\end{align}
The proof is then completed by combining \eqref{eq20181024_1} and \eqref{eq20161222_10}.

\end{proof}
\section{Numerical Tests}\label{sec-6}
In this section, we present two numerical examples to validate the theoretical results. Continuous and discontinuous coefficients are considered in these two problems, respectively. The rates of convergence in the probability space and the physical space are both examined in each test. 
In all the numerical tests, leap-frog time integration is used to achieve energy conservation.

\vspace{.1in}
\noindent{\bf{Test 1 (Continuous coefficient).}} Consider the following wave equation 
\begin{align}
\label{eq-test1}
\frac{\partial^2 u(t,\x,\y)}{\partial t^2} =  \mathrm{div}(a^2(\x,\y)\nabla u(t,\x,\y))\qquad \text{in}\ \mathcal{T}\times\mathcal{D}\times\mathbb{R}^2,
\end{align}
where $\mathcal{T}=[0,T]$ is the time domain, $\D=[0,2]\times[0,2]$ is the physical domain and $\mathbb{R}^2=[-1,1]\times[-1,1]$ is the domain for $\y$. For simplicity, we impose the exact solution (see below) as its boundary conditions.
The coefficient $a$ is defined by
\begin{align*}
a^2(\x,\y)=\frac{2}{(1+\delta y_1)^2+(1+\delta y _2)^2}, 
\end{align*}
where $y_1$ and $y_2$ are two independent random variables with uniform distributions on $[-1,1]$, and $\delta$ is a small number representing the magnitude of perturbation. 
The exact solution is 
\[ u(t,\x,\y)=\cos(\sqrt{2}\pi t)\sin(\pi(1+\delta y_1)x_1)\sin(\pi(1+\delta y_2)x_2). \]

The errors of the numerical solution are defined as:
\begin{align}\label{eq20170302_2}
\|e_u\|_{L^\infty(L^2)} &:=\max_{t\in[0,T]}\bigg(\int_{\D}\E[(u_h-u)^2]d\x\bigg)^{\frac12},\\
\|e_\q\|_{L^\infty(L^2)} &:=\max_{t\in[0,T]}\bigg(\int_{\D}\E[(\mathbf{q}_h-\mathbf{q})^2]d\x\bigg)^{\frac12}.
\end{align}
For simplicity, above we use $L^p(L^q)$ to denote $L^p(\mathcal{T};(L^q(\mathcal{D}))$, where $1\le p, q\le\infty$. 
Table \ref{tab_add_1} shows the $L^\infty(L^2)$ errors and the convergence rates for $u$, $u_x$ and $u_y$, when linear elements are used in LDG discretization. We take $M=15$ ($P=4$) in the gPC expansion, $\delta=0.01$, time step $\Delta t=1.5625\times 10^{-5}$ and final time $T=1.5625\times 10^{-3}$. Second order accuracy can be observed, as expected.
As cubic elements are used in the LDG method, a clear 4-th order can be obtained, as shown in Table \ref{tab_add_2}. 

\begin{table}[htb]
\footnotesize
\begin{center}
\begin{tabular}{|*{22}{c|}}\hline
 &  \multicolumn{2}{|c|}{$u$} & \multicolumn{2}{|c|}{$u_x$}  &  \multicolumn{2}{|c|}{$u_y$}  \\ \hline
$h$
& error & order&  error&  order&  error&  order\\ \hline
$0.5$ & 1.0113E-01  & &2.6454E-01  & &2.6454E-01&\\ \hline
$0.25$ &2.6248E-02  &1.9459 &6.8421E-02  &1.9510 &6.8421E-02 &1.9510 \\ \hline
$0.125$ & 6.6183E-03 &1.9877& 1.7243E-02 &1.9884 &1.7243E-02 &1.9884\\ \hline
$0.0625$ & 1.6580E-03 & 1.9970&4.3192E-03  & 1.9972&4.3192E-03 &1.9972\\ \hline
\end{tabular}
\smallskip
\caption{$L^\infty(L^2)$ errors and order of accuracy with linear elements in LDG method.  $M = 15, \delta=0.01, \Delta t=1.5625\times 10^{-5}, T=1.5625\times 10^{-3}.$} 
\label{tab_add_1} 
\end{center}
\end{table}

\begin{figure}[H]
\centering
\includegraphics[height=2.0in,width=3.0in]{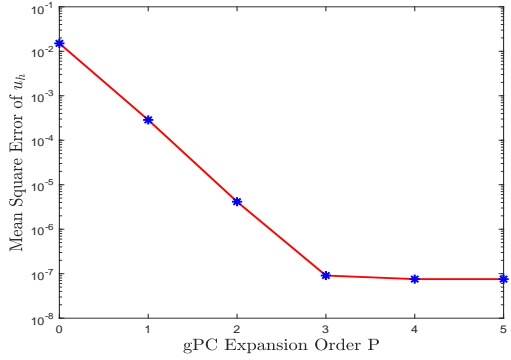}
\caption{$L^\infty(L^2)$ error of $u_h$ with different orders of the gPC expansion. Cubic elements are used in LDG method, with $\delta=0.01, \Delta t=1.5625\times 10^{-5}, T=1.5625\times 10^{-3}.$        
}\label{fig1}
\end{figure}

\begin{table}[htb]
\footnotesize
\begin{center}
\begin{tabular}{|*{22}{c|}}\hline
 &  \multicolumn{2}{|c|}{$u$} & \multicolumn{2}{|c|}{$u_x$}  &  \multicolumn{2}{|c|}{$u_y$}  \\ \hline
$h$
& error & order&  error&  order&  error&  order\\ \hline
$0.5$ &1.2556E-03   & &3.3474E-03  & &3.3474E-03&\\ \hline
$0.25$ &8.0147E-05  &3.9696 &2.1351E-04  &3.9707 &2.1351E-04 &3.9707 \\ \hline
$0.125$ & 5.0356E-06 &3.9924&1.3414E-05  &3.9925 &1.3414E-05 &3.9925\\ \hline
$0.0625$ & 3.1514E-07 &3.9981 &8.4178E-07  &3.9942 &8.4178E-07 &3.9942\\ \hline
\end{tabular}
\smallskip
\caption{$L^\infty(L^2)$ errors and order of accuracy with cubic elements in LDG method.  $M = 15, \delta=0.001, \Delta t=1.5625\times 10^{-5}, T=1.5625\times 10^{-3}.$} 
\label{tab_add_2} 
\end{center}
\end{table}
To test the convergence of gPC expansion in the probability space, we use different orders in the expansion, while fixing the LDG discretization with cubic elements. In Figure \ref{fig1} we observe that the $L^\infty(L^2)$ error decreases exponentially when the order of expansion is increased. 
However, the error saturates for an order larger than 3 because the error from spatial discretization dominates. 

\begin{figure}[H]
\centering
\includegraphics[height=2.2in,width=3.0in]{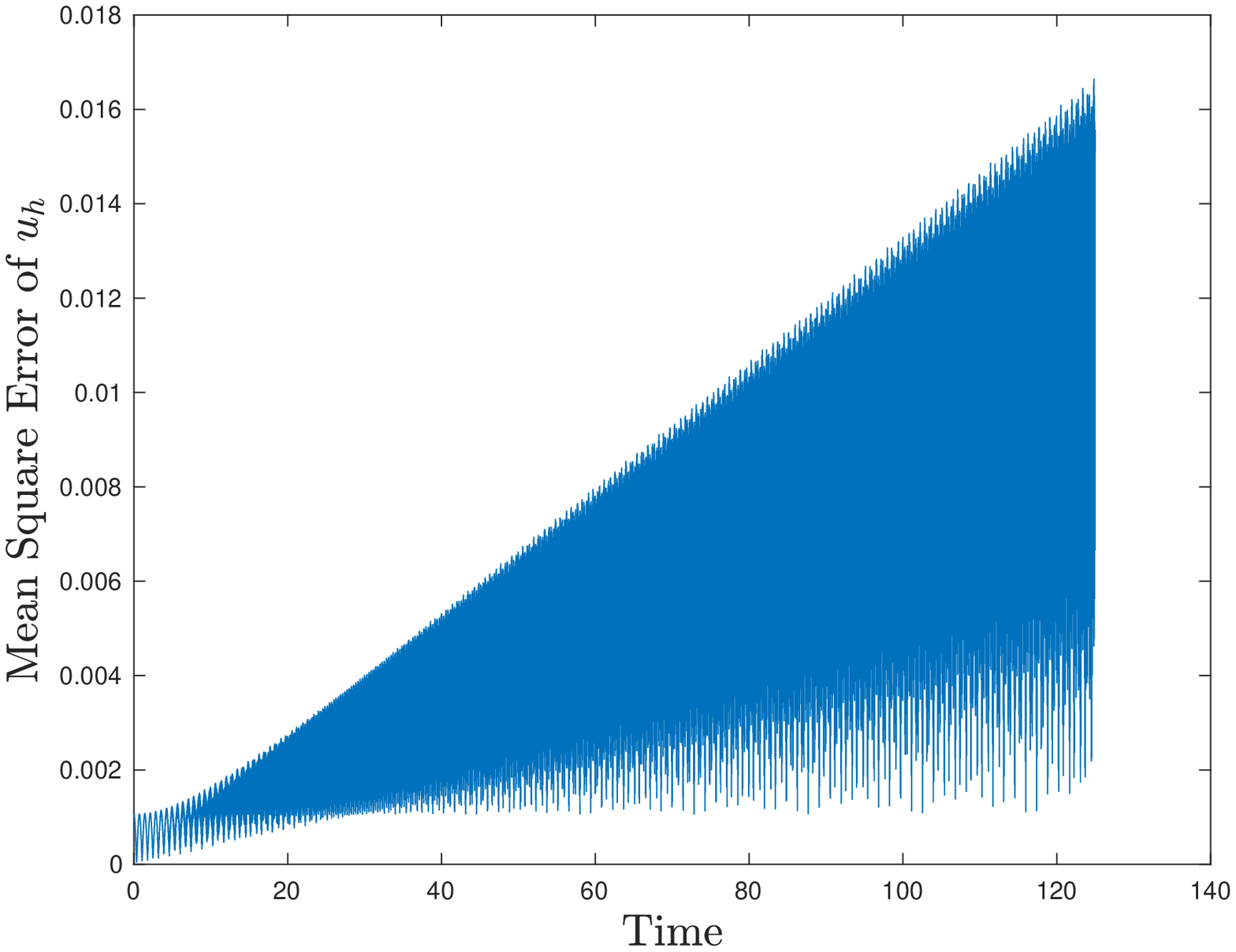}
\includegraphics[height=2.2in,width=3.0in]{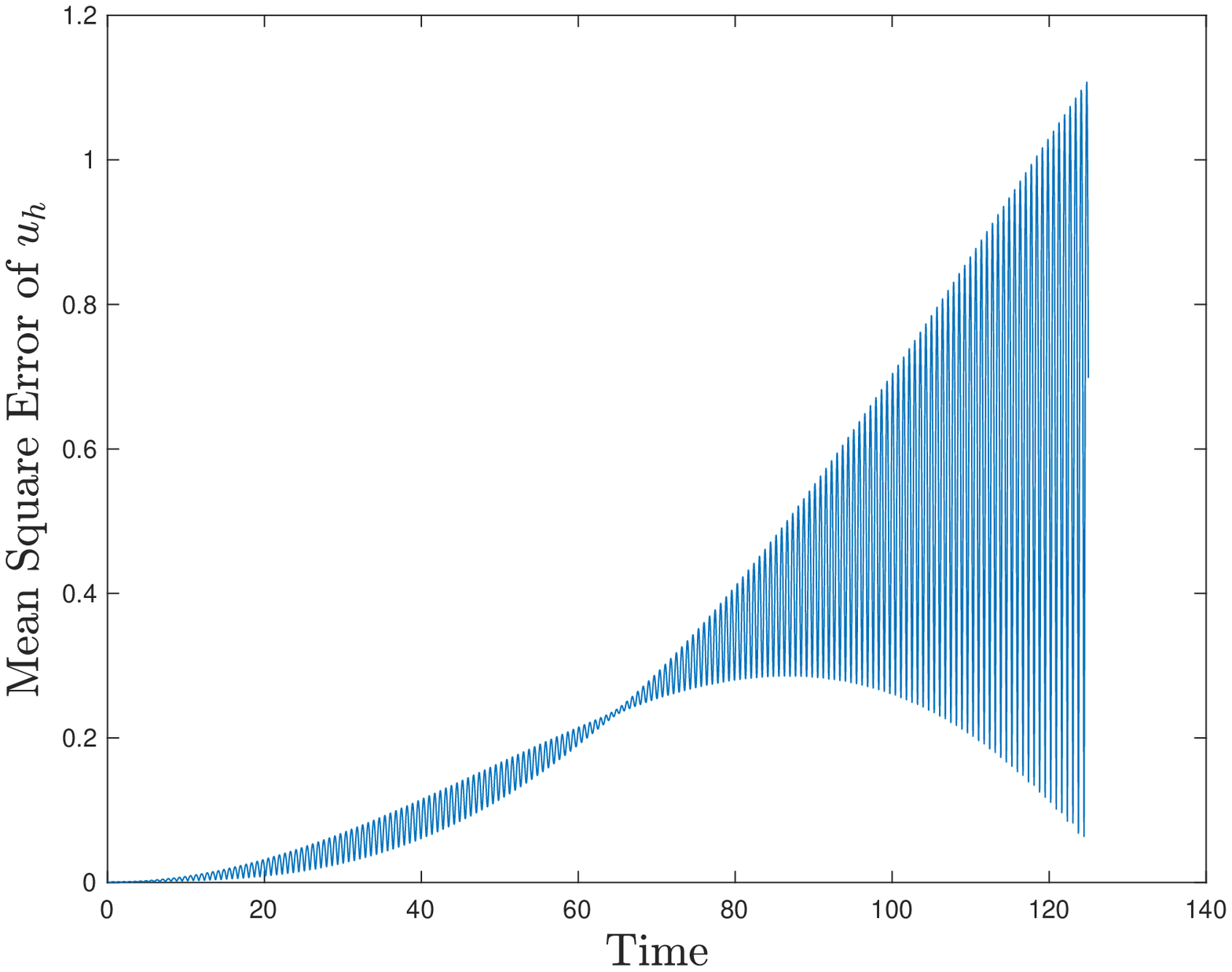}
\caption{Long time $L^\infty(L^2)$ errors of $u_h$ with linear elements in LDG method. Noise with magnitude $\delta=10^{-6}$ is used in the left figure, and  $\delta=10^{-2}$ is used in the right figure. $M=3,\ \Delta t=6.25\times 10^{-5},\  T=125.$}
\label{fig3}
\end{figure}

Next, we demonstrate the advantage of energy conservation property by tracking the errors for a long time simulation.
Figure \ref{fig3} shows the $L^\infty(L^2)$ errors when linear elements are used in LDG and $M=3$ ($P=1$) in gPC expansions. 
In these test cases, both small and large magnitudes of noise ($\delta$) are considered; the time step is $\Delta t=6.25\times10^{-5}$ and the final time is $T=125$. It can be seen that the growth of errors is on average linear or linearly bounded for both cases.

\vspace{.1in}
\noindent{\bf{Test 2 (Discontinuous coefficient).}}
Consider the same equation \eqref{eq-test1} as in Test 1. 
The spatial domain $\D=\D_1\cup\D_2=[-1,1]\times[-1,1]$, with $\D_1=[-1,0]\times[-1,1]$, $\D_2=(0,1]\times[-1,1]$.
The coefficient $a$ is defined by
\begin{equation*}
 a^2(\x,\y)=\begin{cases}
\frac{1}{(1+\delta y_1)^2+(1+\delta y_2)^2} &  \text{in}\ \D_1,\\
\frac{9}{25(1+\delta y_1)^2+9(1+\delta y_2)^2} &  \text{in}\ \D_2,\\
 \end{cases}
 \end{equation*}
 where $y_1$ and $y_2$ are two independent random variables with uniform distributions on $[-1,1]$, and $\delta$ is the magnitude of the noise. 
 We again impose the following exact solution on the boundaries.

\begin{figure}[H]
\centering
\includegraphics[height=2.0in,width=3.0in]{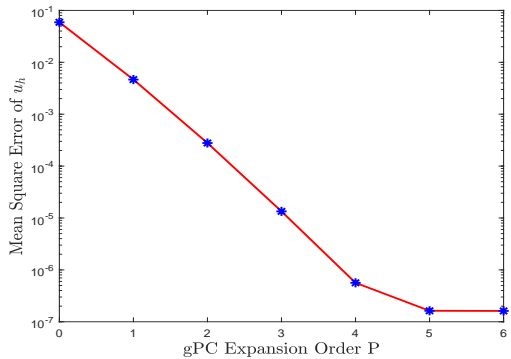}
\caption{$L^\infty(L^2)$ error of $u_h$ with different orders of the gPC expansion. Cubic elements are used in LDG method, with $\delta=0.01, \Delta t=2.5\times 10^{-8}, T=2.5\times 10^{-6}.$ }\label{fig2}
\end{figure}

 The exact solution is 
\begin{equation*}
u(t,\x,\y)=\begin{cases}
\cos(3\pi t)\sin(3\pi(1+\delta y_1)x_1)\sin(3\pi(1+\delta y_2)x_2) &  \text{in}\ \D_1,\\
\cos(3\pi t)\sin(5\pi(1+\delta y_1)x_1)\sin(3\pi(1+\delta y_2)x_2) &  \text{in}\ \D_2.\\
\end{cases}
\end{equation*}
Note that the random coefficient is discontinuous along the vertical line $x=0$. 
Table \ref{tab4} shows the rate of convergence of the numerical method in $L^\infty(L^2)$ norm. We can see that for $u$, $u_x$ and $u_y$ all the errors converge in second order, as expected. In this accuracy test we use $M=15$ ($P=4$) in the gPC expansion with $\delta=0.01$, time step $\Delta t=1.5625\times 10^{-5}$ and final time $T=1.5625\times 10^{-3}$.
Optimal convergence rates are also observed for high order cubic elements, as shown in Table \ref{tab6}. 
In this test, $\delta=0.001$, $\Delta t=2.5\times 10^{-8}$ and $T=2.5\times 10^{-6}$ are used.
In Figure \ref{fig2}, we show that given a fixed spatial discretization in LDG (with cubic elements), the error in $u$ decreases exponentially as the order of gPC expansion becomes higher and saturates when the spatial error dominates. 

\begin{table}[htb]
\footnotesize
\begin{center}
\begin{tabular}{|*{22}{c|}}\hline
 &  \multicolumn{2}{|c|}{$u$} & \multicolumn{2}{|c|}{$u_x$}  &  \multicolumn{2}{|c|}{$u_y$}  \\ \hline
$h$
& error & order&  error&  order&  error&  order\\ \hline
$0.5$ & 5.3285E-01 & & 3.5918E+00 & &2.9336E+00&\\ \hline
$0.25$ & 3.0264E-01 &0.8161 & 1.7426E+00 &1.0435&1.4927E+00&0.9747 \\ \hline
$0.125$ & 9.2197E-02 &1.7148& 5.0773E-01 &1.7791&4.4723E-01&1.7388\\ \hline
$0.0625$ & 2.4080E-02 &1.9369 & 1.3518E-01 &1.9092&1.1663E-01&1.9391\\ \hline
\end{tabular}
\smallskip
\caption{$L^\infty(L^2)$ errors and order of accuracy with linear elemtnes in LDG method.  $M = 15, \delta=0.01, \Delta t=1.5625\times 10^{-5}, T=1.5625\times 10^{-3}.$} 
\label{tab4} 
\end{center}
\end{table}
\begin{table}[htb]
\footnotesize
\begin{center}
\begin{tabular}{|*{22}{c|}}\hline
 &  \multicolumn{2}{|c|}{$u$} & \multicolumn{2}{|c|}{$u_x$}  &  \multicolumn{2}{|c|}{$u_y$}  \\ \hline
$h$
& error & order&  error&  order&  error&  order\\ \hline
$0.5$ & 2.0522E-01 & & 1.2190E+00 & &1.0025E+00&\\ \hline
$0.25$ & 2.1785E-02 &3.2358 & 1.2041E-01 &3.3397&1.0582E-01&3.2439 \\ \hline
$0.125$ & 1.5483E-03 &3.8146& 8.4204E-03 &3.8379&7.5177E-03&3.8152\\ \hline
$0.0625$ & 9.9907E-05 &3.9540 & 5.4516E-04 &3.9491&4.8506E-04&3.9541\\ \hline
\end{tabular}
\smallskip
\caption{$L^\infty(L^2)$ errors and order of accuracy with cubic elements in LDG method. $M = 15, \delta=0.001,  \Delta t=2.5\times 10^{-8}, T=2.5\times 10^{-6}.$} 
\label{tab6} 
\end{center}
\end{table}

\begin{figure}[th]
\centering
\includegraphics[height=2.2in,width=3.0in]{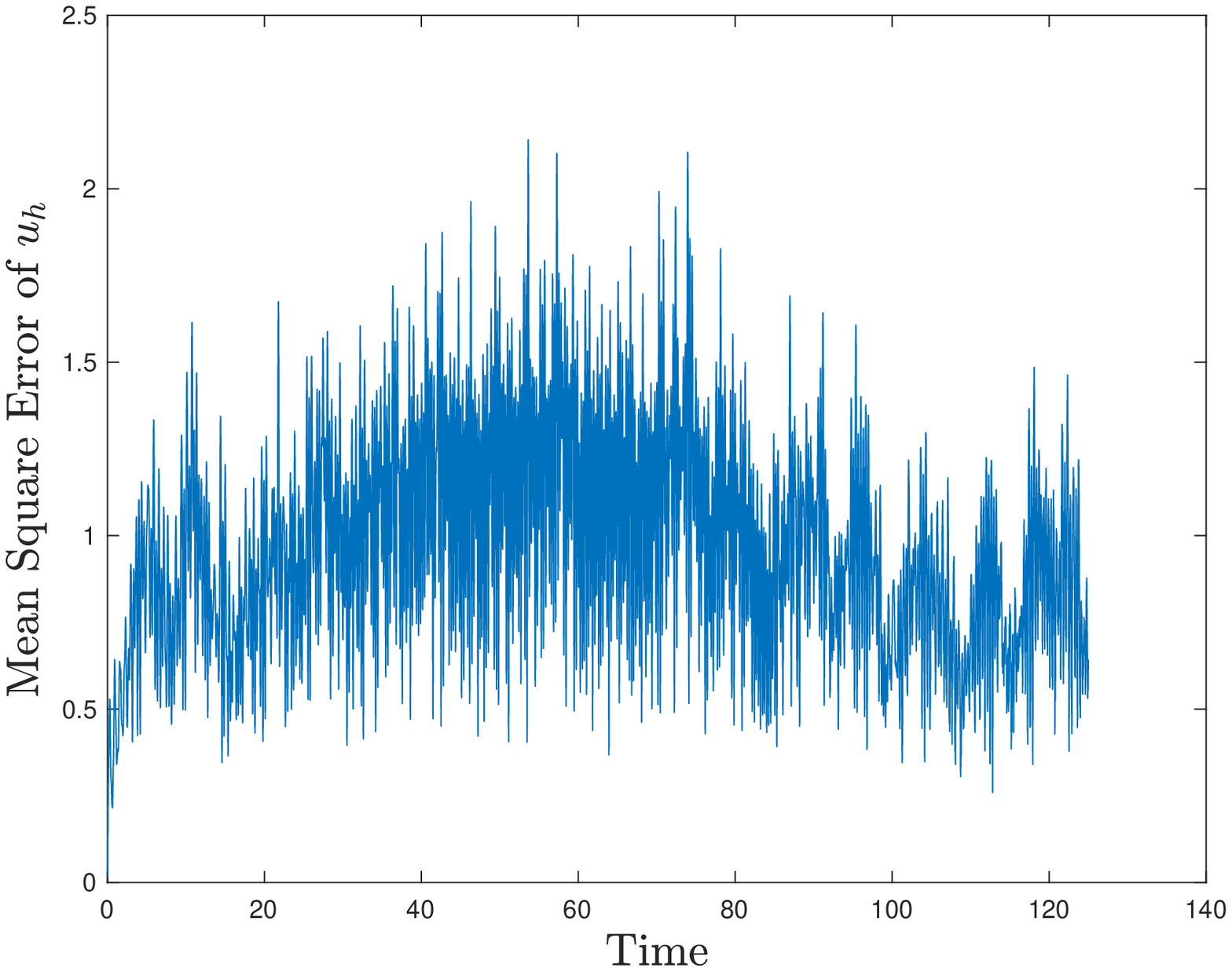}
\includegraphics[height=2.2in,width=3.0in]{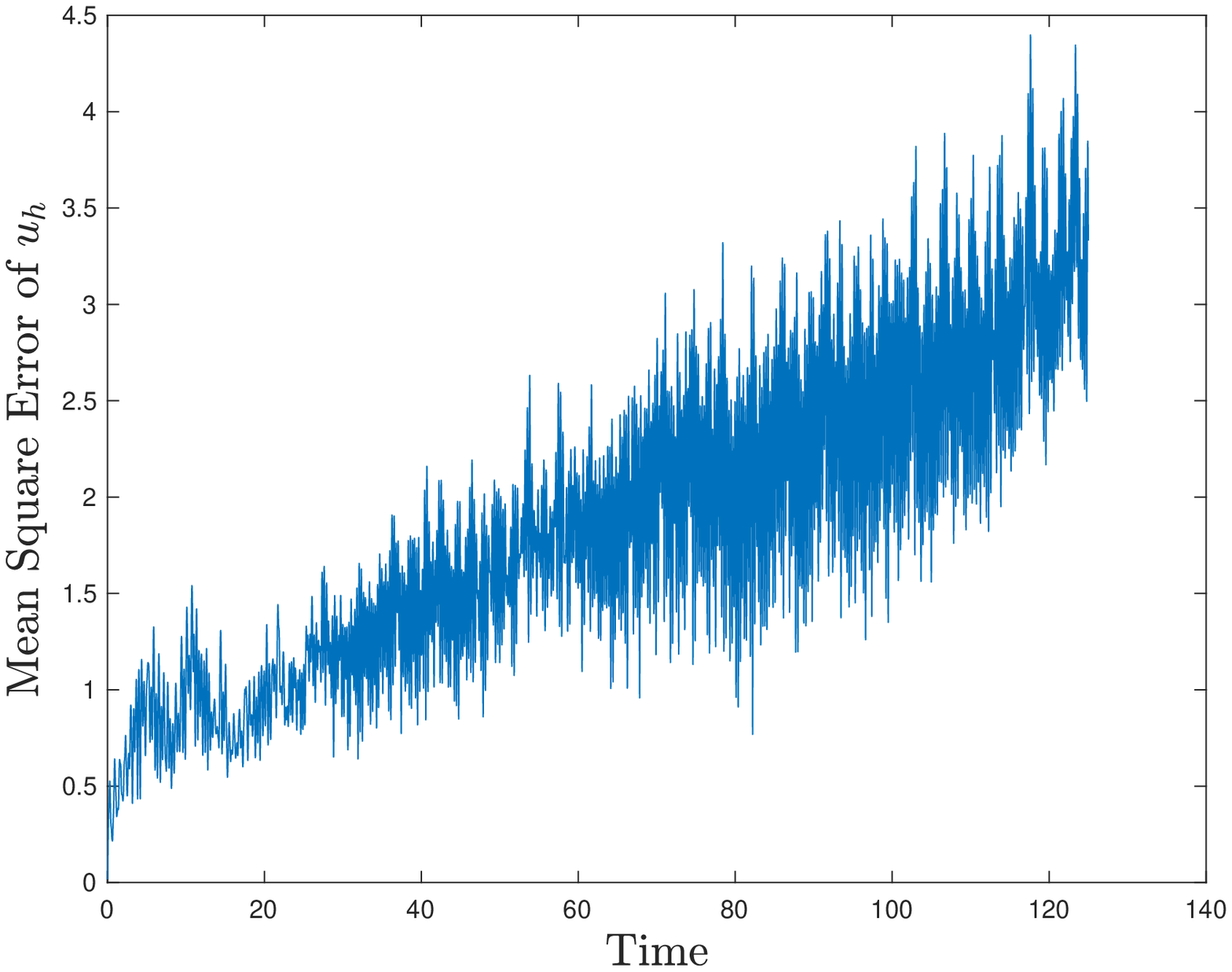}
\caption{Long time $L^\infty(L^2)$ errors of $u_h$ with $P^1$ in LDG method. Smaller noise $\delta=10^{-6}$ is used in the left graph, and bigger noise $\delta=10^{-2}$ is used in the right graph. $M=3,\ \Delta t=6.25\times 10^{-5},\  T=125.$}
\label{fig4}
\end{figure}

Figure \ref{fig4} shows the $L^\infty(L^2)$ errors when linear elements are used in LDG with $M=3$ ($P=1$) in gPC expansions. In these test cases, we consider $\delta=10^{-6}$ and $10^{-2}$, with the time step being $\Delta t=6.25\times10^{-5}$ and final time is $T=125$. The errors appear to be large because we used $M=3$ to save the computational time; however, the errors for both large and small $\delta$'s are linearly bounded as expected from the theoretical results.

\section{Concluding Remarks}

In this paper, we have presented a numerical scheme for solving second-order wave equation with random wave speed coefficient.
Our method is based on gPC expansion with stochastic Galerkin method for probability space, and LDG discretization for physical space.
We are able to show the energy conserving property of the proposed method in both semi-discrete form and  fully-discrete form when 
leap-frog time discretization is used.
The error estimate shows that the convergence of the scheme is optimal, and the grow of the error is at most linear in time.
Taken together, the numerical solution will benefit from these properties and have small shape error (including both dissipative and dispersive errors) and phase error after long time 
integration. Our numerical tests further validate the theoretical findings.



\end{document}